\documentclass[12pt]{article}   %  for old version
\usepackage{latexsym}
\usepackage{amsbsy}
\usepackage{amsmath}   %  for  align  environment used by KM
\usepackage{graphicx}
\usepackage{graphics}
\usepackage{tikz}
\usepackage{pgf}
\usepackage{pgflibraryarrows}
\usepackage{enumerate}   % http://ctan.org/pkg/enumerate
\usepackage{stmaryrd}   %  for blackboard-bold-style square brackets  [[  ...  ]]

\usepackage{amssymb}

\PassOptionsToPackage{hyphens}{url}
\usepackage{hyperref}

\usetikzlibrary{calc}   %  for co-ord calculations

\newtheorem{thm}{Theorem}

\newtheorem{propn}[thm]{Proposition}

\newcommand{\pf}{\noindent\textit{Proof.\ \ }}
\newcommand{\eopf}{\hfill\hspace*{4pt}\hfill
\fbox{\rule[-1pt]{0pt}{4pt}\hspace*{2pt}}}
\definecolor{greyish}{gray}{0.7}
\definecolor{greenish}{rgb}{0,0.5,0}
\definecolor{lightblue}{rgb}{0.5,0.5,1}
\definecolor{yellowish}{rgb}{1,1,0}
\definecolor{oldgold}{rgb}{0.85,.66,0}
\definecolor{gold}{rgb}{1,0.843,0}
\definecolor{mediumgold}{rgb}{0.75,0.632,0}
\definecolor{darkgold}{rgb}{0.5,0.4215,0}
\definecolor{ochre}{rgb}{0.80,0.47,0.13}
\definecolor{darkbrown}{rgb}{0.5,0.25,0}

\newcommand{\contract}{/\hspace{-3pt}/}
\newcommand{\delete}{\setminus\hspace{-5pt}\setminus}
\newcommand{\minor}[1]{\,\|_{{}_{#1}}}
\newcommand{\loopy}{\hbox{\rm loop}}
\newcommand{\coloopy}{\hbox{\rm coloop}}
\newcommand{\scriptloopy}{\hbox{\scriptsize\rm loop}}
\newcommand{\scriptcoloopy}{\hbox{\scriptsize\rm coloop}}
\newcommand{\dom}{\mathop{\mathrm{dom}}}

\newcommand{\Go}{\hbox{\rm Go}}

\title{Graph polynomials: some questions on the edge}
\author{Graham Farr   \\
Department of Data Science and AI   \\
Faculty of I.T.   \\
Monash University   \\
Australia   \\
\href{mailto:Graham.Farr@monash.edu}{\texttt{Graham.Farr@monash.edu}}   \\
\hspace*{1cm}   \\
\tiny
\and
Kerri Morgan   \\
School of Science (Mathematical Sciences)   \\
RMIT University   \\
Australia   \\
\href{mailto:Kerri.Morgan@rmit.edu.au}{\texttt{Kerri.Morgan@rmit.edu.au}}   \\
\hspace*{1cm}   \\
\tiny
}

\date{22 June 2024}

\begin{document}

\maketitle

\begin{abstract}
%~   %  seems to be needed for published version?
We raise some questions about graph polynomials, highlighting concepts and
phenomena that may merit consideration in the development of a general theory.
Our questions are mainly of three types:
When do graph polynomials have reduction relations
(simple linear recursions based on local operations), perhaps in a wider class of combinatorial
objects?
How many levels of reduction relations does a graph polynomial need
in order to express it in terms of trivial base cases?
For a graph polynomial,
how are properties such as equivalence and factorisation
reflected in the structure of a graph?
We illustrate our discussion with a variety of graph polynomials and other invariants.
This leads us to reflect on the historical origins of graph polynomials.
We also introduce some new polynomials based on partial colourings of graphs and
establish some of their basic properties.
\end{abstract}

\section{Introduction}

J\'anos Makowsky started his research career in mathematical logic over half a century ago.
For the last two decades, he has brought many concepts and results from that
field --- along with the perspective it offers --- to bear on the study of graph polynomials.
This has led to significant new theorems, links between topics, fresh viewpoints,
and deeper understanding.

His contributions to graph polynomials, with a wide set of collaborators,
include: complexity classifications for specific computational problems \cite{blaeser-dell-makowsky2010,lotz-makowsky2004,makowsky2005,makowsky-marino2003}
(see also \cite{kotek-makowsky2022});
refining the models of computation used to study the computational complexity of problems
concerning graph polynomials \cite{makowsky-meer2000,kotek-makowsky-ravve2013};
introducing formal logical definitions of graph polynomials,
especially using second-order logic (SOL) and variants thereof, and using concepts and tools
from logic to develop their theory
at a very general level \cite{courcelle-makowsky-rotics01,makowsky04,makowsky06,kotek-makowsky-zilber08,godlin-katz-makowsky2012};
developing a general framework for studying the distinguishing power of graph polynomials \cite{makowsky-ravve-blanchard2014,makowsky-ravve2016,kotek-makowsky-ravve2018};
and showing that the locations of zeros of a graph polynomial is ``not a semantic property'',
in that it derives more from the algebraic form of the polynomial than from the way it partitions
the set of all graphs into
equivalence classes \cite{makowsky-ravve-blanchard2014,makowsky-ravve-kotek2019}.
Some of his own reflections on these topics may be found in \cite{makowsky-ravve-kotek2019,makowsky2023,makowsky2024}.

A pervasive characteristic of his work has been to put specific graph polynomials in
a wider mathematical context --- to see the forest as well as the trees, to quote his own quotation
of Einstein \cite{makowsky2012}.
One manifestation of this has been his 
advocacy for the development of a general theory of graph polynomials
and his own work in that direction \cite{makowsky06,makowsky08,makowsky-ravve-kotek2019}.
He was a co-organiser of the Dagstuhl Seminar on ``Graph Polynomials: Towards a Comparative Theory''
in 2016 \cite{ellis-monaghan-etal2016}.
At the Dagstuhl Seminar on ``Comparative Theory for Graph Polynomials''
in 2019 \cite{ellis-monaghan-etal2019}, he helped lead an informal group working on the distinguishing
power of graph polynomials.

In this paper we propose some polynomials,
topics and questions that may merit further consideration in
the development of a general theory of graph polynomials.  To do this, we journey to the edge
of the territory covered by the current theory, meeting some polynomials that seem to lie near,
or beyond, that frontier, as well as some that are very familiar and are well covered by the theory.

We begin with some general definitions and notation in \S\ref{sec:defn-notn}
and set the scene in \S\ref{sec:some-gph-polys} by defining all the polynomials we will discuss.
This includes the introduction of some new graph polynomials related to partial colourings.
In \S\ref{sec:origins} we reflect on the origins of graph polynomials.
Some of the polynomials we introduce
are then used to illustrate the questions raised in the next four sections.
In \S\ref{sec:redn-relns} we consider the widespread phenomenon of graph polynomials having
reduction relations (i.e., simple
recursive relations based on local modifications), pointing out that even graph polynomials
that do not seem
to have such a relation will often be found to have one within a wider class of objects.
In \S\ref{sec:levels} we discuss a notion of ``levels'' in these reduction relations,
where a graph polynomial
can be reduced to a large sum of polynomials of reduced objects of some kind, which in turn may
each be reduced to another large sum using another relation, and so on.
In \S\ref{sec:gph-polys-qn} we discuss the relationship between the algebraic properties of a graph polynomial
and the structure of the graph.
In \S\ref{sec:certs} we pose questions about certificates, a tool for studying equivalence and
factorisation of graph polynomials.

Some of the material in this paper was presented by the first author in a talk of the same
title at the Dagstuhl Seminar on ``Graph Polynomials: Towards a Comparative Theory''
in June 2016 \cite{ellis-monaghan-etal2016}.

\section{Definitions and notation}
\label{sec:defn-notn}

Throughout, $G=(V,E)$ is a graph with $n$ vertices and $m$ edges.
The number of components of $G$ is denoted by $k(G)$.
If $X\subseteq E$ then $V(X)$ denotes the set of vertices of $G$ that are incident with at least
one edge in $X$ (overloading the $V(\_)$ notation slightly).
We write $\nu(X)$ for the number of vertices of $G$ that meet an edge in $X$
(succinctly: $\nu(X)=|V(X)|$).
The number of components of $(V,X)$ is denoted by $k(X)$, while
the number of components of $(V(X),X)$ is denoted by $\kappa(X)$.  The former count
includes isolated vertices, while the latter count excludes them: $k(X) = \kappa(X)+n-\nu(X)$.
We write $\rho(X)$ for the \textit{rank} of $X$,
given by $\rho(X)=\nu(X)-\kappa(X)=|V|-k(X)$, and $\rho(G):=\rho(E(G))$.
For any function $r:2^E\rightarrow\mathbb{R}$, its \textit{dual} $r^*$ is defined by
$r^*(X)=|X|+r(E\setminus X)-r(E)+r(\emptyset)$.  When $\rho$ is the rank function of $G$
(and therefore the rank function of its cycle matroid),
$\rho^*$ is the rank function of the dual of the cycle matroid of $G$.

If $U\subseteq V$ then $G[U]$ is the subgraph of $G$ induced by $U$.

The disjoint union of two graphs $G$ and $H$ is denoted by $G\sqcup H$.

If $e\in E$ then $G\setminus e=(V,E\setminus\{e\})$ is the graph obtained from $G$
by deleting edge $e$ and $G/e$ is the graph obtained by contracting edge $e$,
i.e., deleting $e$ and then identifying its former endpoints.  If $u,v\in V$ with $uv\not\in E$ then 
$G+uv=(V,E\cup\{uv\})$ is the graph obtained by adding an edge between $u$ and $v$ in $G$ and $G/uv$ is obtained from $G$ by identifying vertices $u$ and $v.$

A \textit{coloop} of $G$ is an edge $e$ such that $k(G\setminus e)>k(G)$.
This is often called a \textit{bridge} and sometimes an \textit{isthmus}.   %  This sentence new in v3.

The maximum degree of $G$ is denoted by $\Delta(G)$.

A \textit{null graph} is a graph with no edges.

%Let $H_1$ and $H_2$ be graphs, with $u_i,v_i\in V(H_i)$ ($i=1,2$).
%Suppose that graph $G$ can be formed from $H_1$ and $H_2$
%by identifying $u_1$ with $v_1$ and $u_2$ with $v_2$,
%and that $G^{\sim}$ can be formed from the same two graphs but by identifying
%$u_1$ with $v_2$ and $u_2$ with $v_1$.  Then we say that $G^{\sim}$ is formed from
%$G$ by a \textit{Whitney twist}.  Graphs $G$ and $G'$ are \textit{2-isomorphic}
%if a graph isomorphic to $G'$ can be constructed from $G$ by a sequence of Whitney
%twists.
 
A \textit{graph invariant} is a function defined on all graphs that depends only on
the isomorphism class of the graph.

We write $[k] = \{1,2, \ldots, k\}$.

The falling factorial $x(x-1)\cdots(x-k+1)$ is denoted by $(x)_k$.

Let $\Lambda$ be a set whose members we will call \textit{colours},
and let $\lambda\in\mathbb{N}$.
A $\Lambda$\textit{-assignment} of $G$ is a function $f:V\rightarrow\Lambda$.
A $\lambda$\textit{-assignment} is a $[\lambda]$-assignment.
A \textit{partial $\Lambda$-assignment} is a function $f:W\rightarrow\Lambda$ where
$W\subseteq V$.
The vertices of $W$ and $V\setminus W$ are \textit{coloured} and \textit{uncoloured}, respectively,
by $f$.
A \textit{partial $\lambda$-assignment} is a partial $[\lambda]$-assignment.
For every $k\in\Lambda$, its \textit{colour class} $C(k)=C_f(k)$
under a partial $\Lambda$-assignment $f$ is given by
$C(k)=f^{-1}(k)=\{v\in V:f(v)=k\}$.
Every colour class $C(k)$ induces a subgraph $G[C(k)]$ of $G$.
A \textit{chromon}, or \textit{monochromatic component}, of $(G,f)$ is a component of
$G[C(k)]$ for some $k\in[\lambda]$.

Every partial $\lambda$-assignment $f$ is determined by its $\lambda$-tuple
$(C_f(i))_{i=1}^{\lambda}$ of colour classes:
given a $\lambda$-tuple $(C_i)_{i=1}^{\lambda}$ of mutually
disjoint subsets of $V$, we can define a partial $\lambda$-assignment $f:\bigcup_{i=1}^{\lambda}C_i\rightarrow[\lambda]$ by $f(v)=i$ for each $v\in C_i$ and $i\in[\lambda]$; this then satisfies
$C_f(i)=C_i$ for all $i\in[\lambda]$.

A colour class is \textit{proper} if it is a stable set in $G$, i.e., no two of its vertices are adjacent in $G$.
A partial $\lambda$-assignment is a \textit{partial $\lambda$-colouring} if every colour class is proper;
this is regardless of whether or not the partial $\lambda$-assignment can be
extended to a $\lambda$-colouring of $G$.

An \textit{extension} of a partial $\lambda$-assignment $f$
is a partial $\lambda$-assignment $g$ such that $f\subseteq g$, which is equivalent to
requiring that $\dom f\subseteq\dom g$ and $f(v)=g(v)$ for all $v\in\dom f$.

\section{Some graph polynomials}
\label{sec:some-gph-polys}

Papers by J\'anos often include observations about collections of specific graph polynomials
as motivation for studying more general phenomena.
In a similar spirit, we now discuss an eclectic collection of graph invariants ---
some old, some new; mostly polynomials, some not --- to help motivate
some of the questions we ask.
The reader can skip those sections treating graph polynomials with which they are familiar.

More comprehensive collections of graph polynomials may be found in \cite{ellis-monaghan-moffatt2022,tittmann2024}.

\subsection{Tutte-Whitney polynomials}

The preeminent graph polynomial is arguably the Tutte polynomial,
due to W.T. Tutte \cite{tutte47,tutte54} and closely related to
the Whitney rank generating function \cite{whitney32d}.

The \textit{Whitney rank generating function} $R(G;x,y)$ of a graph $G$ is the bivariate polynomial
defined by
\[
R(G;x,y)
~=~ \sum_{X\subseteq E} x^{\rho(E)-\rho(X)}y^{\rho^*(E)-\rho^*(E\setminus X)}
~=~ \sum_{X\subseteq E} x^{\rho(E)-\rho(X)}y^{|X|-\rho(X)}
\]

The \textit{Tutte polynomial} $T(G;x,y)$ may be defined by
\begin{equation}
\label{eq:tgxy}
T(G;x,y) ~=~ \left\{
\begin{array}{ll}
1,  &  \hbox{if $G$ has no edges;}   \\
x\,T(G\setminus e;x,y),  &  \hbox{if $e$ is a coloop in $G$;}   \\
y\,T(G/e;x,y),  &  \hbox{if $e$ is a loop in $G$;}   \\
T(G\setminus e;x,y) + T(G/e;x,y),~~~~  &  \hbox{otherwise,}   \\
\end{array}
\right.
\end{equation}
for any edge $e\in E(G)$.
The polynomials of Whitney and Tutte are related by a simple coordinate translation:
$T(G;x,y)=R(G;x-1,y-1)$ \cite{tutte47,tutte54}.

The recurrence (\ref{eq:tgxy}) is the most fundamental property of the
Tutte polynomial, and one of Tutte's major conceptual contributions
in \cite{tutte47} was to base the theory on this relation, a fundamental
conceptual advance upon the pioneering work of Whitney \cite{whitney32d}.
For a history of these polynomials, see \cite{farr2022}.

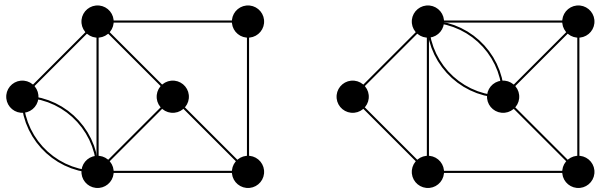
\begin{figure}
\begin{center}
\begin{tikzpicture}
\path (0,1) node(a) {} ;
\path (1,0) node(b) {} ;
\path (1,2) node(c) {} ;
\path (2,1) node(d) {} ;
\path (3,0) node(e) {} ;
\path (3,2) node(f) {} ;
\draw[thick,fill=black] (a) circle (0.2cm) ;
\draw[thick,fill=black] (b) circle (0.2cm) ;
\draw[thick,fill=black] (c) circle (0.2cm) ;
\draw[thick,fill=black] (d) circle (0.2cm) ;
\draw[thick,fill=black] (e) circle (0.2cm) ;
\draw[thick,fill=black] (f) circle (0.2cm) ;
\draw[thick] (a) arc(180:270:1cm) -- (b) ;
\draw[thick] (a) arc(90:0:1cm) -- (b) ;
\draw[thick] (a) -- (c) ;
\draw[thick] (b) -- (c) ;
\draw[thick] (b) -- (d) ;
\draw[thick] (c) -- (d) ;
\draw[thick] (b) -- (e) ;
\draw[thick] (d) -- (e) ;
\draw[thick] (c) -- (f) ;
\draw[thick] (e) -- (f) ;
\end{tikzpicture}
~~~~~
\begin{tikzpicture}
\path (0,1) node(a) {} ;
\path (1,0) node(b) {} ;
\path (1,2) node(c) {} ;
\path (2,1) node(d) {} ;
\path (3,0) node(e) {} ;
\path (3,2) node(f) {} ;
\draw[thick,fill=black] (a) circle (0.2cm) ;
\draw[thick,fill=black] (b) circle (0.2cm) ;
\draw[thick,fill=black] (c) circle (0.2cm) ;
\draw[thick,fill=black] (d) circle (0.2cm) ;
\draw[thick,fill=black] (e) circle (0.2cm) ;
\draw[thick,fill=black] (f) circle (0.2cm) ;
\draw[thick] (a) -- (b) ;
\draw[thick] (a) -- (c) ;
\draw[thick] (b) -- (c) ;
\draw[thick] (d) -- (f) ;
\draw[thick] (c) arc(180:270:1cm) -- (d) ;
\draw[thick] (c) arc(90:0:1cm) -- (d) ;
\draw[thick] (b) -- (e) ;
\draw[thick] (d) -- (e) ;
\draw[thick] (c) -- (f) ;
\draw[thick] (e) -- (f) ;
\end{tikzpicture}
\end{center}
\caption{The two Gray graphs, from \cite{tutte74}.}
\label{fig:gray-graphs}
\end{figure}

A graph $G$ is \textit{Tutte equivalent} to a graph $H$ if $T(G; x,y)=T(H; x,y)$ \cite{tutte74,whitney32d,whitney33}.  Graphs with isomorphic cycle matroids are Tutte equivalent,
since the Tutte polynomial of a graph depends only on its cycle matroid.
But, as an equivalence relation on graphs,
Tutte equivalence is coarser than cycle matroid isomorphism.
Tutte \cite{tutte74} gave two graphs, due to M.C. Gray, which are not
isomorphic, and do not even have isomorphic cycle matroids,
but which have the same Tutte polynomial:
see Figure~\ref{fig:gray-graphs}.  Many other such pairs are known.

The most famous specialisation of the Tutte polynomial is the \textit{chromatic polynomial}
$P(G;q)=(-1)^{\rho(G)}q^{k(G)}T(G;-q+1,0)$ introduced by Birkhoff in \cite{birkhoff1912-1913}.
For $q\in\mathbb{N}$, it gives the number of $q$-colourings of $G$.
Two graphs are \textit{chromatically equivalent} if they have the same chromatic polynomial.

Tutte-Whitney polynomials have been generalised from graphs
to many other combinatorial objects \cite{ellis-monaghan-moffatt2022,farr07a}.
For structures on which deletion and contraction operations exist,
Krajewski, Moffatt and Tanasa \cite{krajewski-moffatt-tanasa2018}
show how to use Hopf algebras
to define a polynomial that may reasonably be called the Tutte polynomial for those structures and
which satisfies a deletion-contraction relation.

\subsection{Some partition functions}
\label{sec:partn-fns}

We introduce the partition functions of three interaction models on graphs: the Ising, Potts
and Ashkin-Teller models.  Our main focus later will be on the Ashkin-Teller model.

The set of edges of a graph $G$ whose endpoints receive the same colour under
a $q$-assignment $f$ is denoted by $E^+(f)$;
these edges are sometimes called \textit{bad} since they are not properly
coloured in the sense of graph colouring.  The set of edges whose endpoints are differently coloured
under $f$ is denoted by $E^-(f)$; these edges are sometimes called \textit{good}.
Note that $E^+(f)\cup E^-(f)=E$.  So the positive and negative signs in superscripts here represent
``bad'' and ``good'', respectively, which is the opposite of their connotations in ordinary English
usage.

With this notation, we may write the chromatic polynomial as
\[
P(G;q) = \sum_{f:V\rightarrow[q]} 0^{|E^+(f)|} 1^{|E^-(f)|}
\]
where $0^k$ is taken to be 1 if $k=0$ and 0 otherwise, so that only proper colourings contribute
to the sum, with all proper colourings counted once.

Suppose now that we do not penalise improper colourings
so drastically, but just weight colourings according to an exponential function of the number of
good edges they have.  This gives the \textit{Potts model}, introduced in \cite{potts52}
and generalising the $q=4$ case introduced by Ashkin and Teller in \cite{ashkin-teller43}.
(See \cite[\S34.12]{farr2022}.
The model is called the Ashkin-Teller-Potts model in \cite[\S4.4]{welsh93}.)
The \textit{Potts model partition function} is given by\footnote{Sometimes, it is defined instead
as $\sum_{f:V\rightarrow[q]}  e^{K \cdot |E^+(f)|}$, which may be written $e^{K|E|}\sum_{f:V\rightarrow[q]}  e^{-K \cdot |E^-(f)|}$.  So the only change is an extra prefactor $e^{K|E|}$.  See, e.g., \cite[\S2]{beaudin-etal2010}.}
\[
Z_{\scriptsize\hbox{Potts}}(G;K,q) = \sum_{f:V\rightarrow[q]}  e^{-K \cdot |E^-(f)|} .
\]
This is a polynomial in $e^{-K}$ and is known to be a partial evaluation of
the Tutte polynomial \cite{essam71,fortuin-kasteleyn72}.  The relationship between them
can be expressed more simply
in terms of the Whitney rank generating function:
\[
Z_{\scriptsize\hbox{Potts}}(G;K,q) = q^{k(G)
}\,(e^K-1)^{\rho(G)}e^{-K|E|} R(G; \frac{q}{e^K-1}, e^K-1) .
\]

The $q=2$ case of the Potts model partition function is mathematically almost the same as
the Ising model partition function \cite{ising1925}.
In the Ising model, colours take values in $\{\pm1\}$, and we sum over all
$\{\pm1\}$-assignments $\sigma:V\rightarrow\{\pm1\}$.
For a given $\sigma:V\rightarrow\{\pm1\}$,
the edge $uv$ belongs to the set $E^{\sigma(u)\sigma(v)}(f)$.
The \textit{Ising model partition function} is
\begin{eqnarray*}
Z_{\scriptsize\hbox{Ising}}(G;K)
& = &
\sum_{\sigma:V\rightarrow\{\pm1\}}  e^{K \sum_{uv\in E}\sigma(u)\sigma(v)}   \\
& = &
\sum_{\sigma:V\rightarrow\{\pm1\}}  e^{K \cdot |E^+(\sigma)|-K \cdot |E^-(\sigma)|}   \\
& = &
e^{K|E|}\sum_{\sigma:V\rightarrow\{\pm1\}}  e^{-2K \cdot |E^-(\sigma)|}   \\
& = &
e^{K|E|} \, Z_{\scriptsize\hbox{Potts}}(G;2K,2) .
\end{eqnarray*}

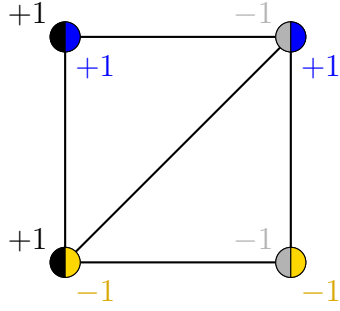
\begin{figure}[h]
\begin{center}
\begin{tikzpicture}
\draw[thick] (0,0) -- (3,0) ;
\draw[thick] (0,0) -- (0,3) ;
\draw[thick] (0,0) -- (3,3) ;
\draw[thick] (0,3) -- (3,3) ;
\draw[thick] (3,0) -- (3,3) ;
\draw[fill=white,opaque] (0,0) circle (0.2cm) ; 
\draw[fill=white,opaque] (3,0) circle (0.2cm) ; 
\draw[fill=white,opaque] (0,3) circle (0.2cm) ; 
\draw[fill=white,opaque] (3,3) circle (0.2cm) ;

\draw (-0.5,0.3) node {$+1$} ;
\draw (0.4,-0.4) node {\textcolor{oldgold}{$-1$}} ;
\draw (2.5,0.3) node {\textcolor{greyish}{$-1$}} ;
\draw (3.4,-0.4) node {\textcolor{oldgold}{$-1$}} ;
\draw (-0.5,3.3) node {$+1$} ;
\draw (0.4,2.6) node {\textcolor{blue}{$+1$}} ;
\draw (2.5,3.3) node {\textcolor{greyish}{$-1$}} ;
\draw (3.4,2.6) node {\textcolor{blue}{$+1$}} ;
\draw[fill=black,opaque] (0,-0.2) arc (270:90:0.2cm) -- cycle ;
\draw[fill=gold,opaque] (0,-0.2) arc (-90:90:0.2cm) -- cycle ;
\draw[fill=greyish,opaque] (3,-0.2) arc (270:90:0.2cm) -- cycle ;
\draw[fill=gold,opaque] (3,-0.2) arc (-90:90:0.2cm) -- cycle ;
\draw[fill=black,opaque] (0,2.8) arc (270:90:0.2cm) -- cycle ;
\draw[fill=blue,opaque] (0,2.8) arc (-90:90:0.2cm) -- cycle ;
\draw[fill=greyish,opaque] (3,2.8) arc (270:90:0.2cm) -- cycle ;
\draw[fill=blue,opaque] (3,2.8) arc (-90:90:0.2cm) -- cycle ;
\end{tikzpicture}
\caption{An Ashkin-Teller model configuration $(\sigma,\tau)$ for a graph.}
\label{fig:ashkin-teller-col}
\end{center}
\end{figure}

The Ashkin-Teller model \cite{ashkin-teller43} extends the Ising model and the
$q=4$ case of the Potts model.
Each vertex $v\in V$ has a pair of colours $\sigma(v)$ and $\tau(v)$, each taking one of the two
values $\pm1$, so the
available colour pairs $(\sigma(v),\tau(v))$ for each vertex are the four pairs $(\pm1,\pm1)$.
We think of the $\{\pm1\}$-assignments
$\sigma:V\rightarrow\{\pm1\}$ and $\tau:V\rightarrow\{\pm1\}$ as two Ising spin functions on $G$.
An example is given in Figure~\ref{fig:ashkin-teller-col}, where the spins at upper left and lower right of each vertex are those assigned by $\sigma$ and $\tau$, respectively.  These spins are also shown as colours on the left and right semicircles in each vertex.

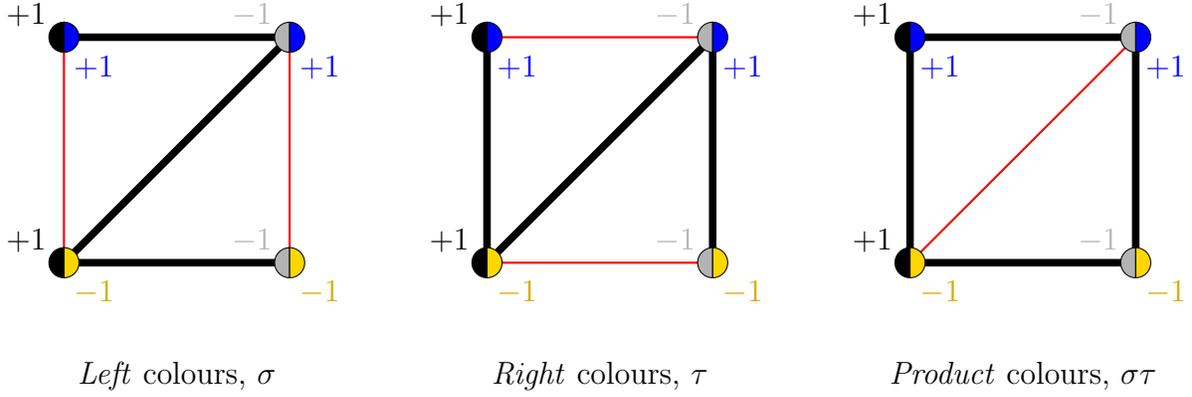
\begin{figure}[t]
\begin{center}
\begin{tikzpicture}
\draw[line width=0.1cm] (0,0) -- (3,0) ;
\draw[line width=0.1cm] (0,0) -- (3,3) ;
\draw[line width=0.1cm] (0,3) -- (3,3) ;
\draw[thick,red] (0,0) -- (0,3) ;
\draw[thick,red] (3,0) -- (3,3) ;
\draw (-0.5,0.3) node {$+1$} ;
\draw (0.4,-0.4) node {\textcolor{oldgold}{$-1$}} ;
\draw (2.5,0.3) node {\textcolor{greyish}{$-1$}} ;
\draw (3.4,-0.4) node {\textcolor{oldgold}{$-1$}} ;
\draw (-0.5,3.3) node {$+1$} ;
\draw (0.4,2.6) node {\textcolor{blue}{$+1$}} ;
\draw (2.5,3.3) node {\textcolor{greyish}{$-1$}} ;
\draw (3.4,2.6) node {\textcolor{blue}{$+1$}} ;
\draw[fill=black,opaque] (0,-0.2) arc (270:90:0.2cm) -- cycle ;
\draw[fill=gold,opaque] (0,-0.2) arc (-90:90:0.2cm) -- cycle ;
\draw[fill=greyish,opaque] (3,-0.2) arc (270:90:0.2cm) -- cycle ;
\draw[fill=gold,opaque] (3,-0.2) arc (-90:90:0.2cm) -- cycle ;
\draw[fill=black,opaque] (0,2.8) arc (270:90:0.2cm) -- cycle ;
\draw[fill=blue,opaque] (0,2.8) arc (-90:90:0.2cm) -- cycle ;
\draw[fill=greyish,opaque] (3,2.8) arc (270:90:0.2cm) -- cycle ;
\draw[fill=blue,opaque] (3,2.8) arc (-90:90:0.2cm) -- cycle ;
\draw (1.5,-1.5) node  {\textit{Left} colours, $\sigma$} ;
\end{tikzpicture}%
~~~~~~~~~~%
\begin{tikzpicture}
\draw[line width=0.1cm] (0,0) -- (0,3) ;
\draw[line width=0.1cm] (0,0) -- (3,3) ;
\draw[line width=0.1cm] (3,0) -- (3,3) ;
\draw[thick,red] (0,0) -- (3,0) ;
\draw[thick,red] (0,3) -- (3,3) ;
\draw (-0.5,0.3) node {$+1$} ;
\draw (0.4,-0.4) node {\textcolor{oldgold}{$-1$}} ;
\draw (2.5,0.3) node {\textcolor{greyish}{$-1$}} ;
\draw (3.4,-0.4) node {\textcolor{oldgold}{$-1$}} ;
\draw (-0.5,3.3) node {$+1$} ;
\draw (0.4,2.6) node {\textcolor{blue}{$+1$}} ;
\draw (2.5,3.3) node {\textcolor{greyish}{$-1$}} ;
\draw (3.4,2.6) node {\textcolor{blue}{$+1$}} ;
\draw[fill=black,opaque] (0,-0.2) arc (270:90:0.2cm) -- cycle ;
\draw[fill=gold,opaque] (0,-0.2) arc (-90:90:0.2cm) -- cycle ;
\draw[fill=greyish,opaque] (3,-0.2) arc (270:90:0.2cm) -- cycle ;
\draw[fill=gold,opaque] (3,-0.2) arc (-90:90:0.2cm) -- cycle ;
\draw[fill=black,opaque] (0,2.8) arc (270:90:0.2cm) -- cycle ;
\draw[fill=blue,opaque] (0,2.8) arc (-90:90:0.2cm) -- cycle ;
\draw[fill=greyish,opaque] (3,2.8) arc (270:90:0.2cm) -- cycle ;
\draw[fill=blue,opaque] (3,2.8) arc (-90:90:0.2cm) -- cycle ;
\draw (1.5,-1.5) node  {\textit{Right} colours, $\tau$} ;
\end{tikzpicture}%
~~~~~~~~~~%
\begin{tikzpicture}
\draw[line width=0.1cm] (0,0) -- (3,0) ;
\draw[line width=0.1cm] (3,0) -- (3,3) ;
\draw[line width=0.1cm] (0,0) -- (0,3) ;
\draw[line width=0.1cm] (0,3) -- (3,3) ;
\draw[thick,red] (0,0) -- (3,3) ;
\draw (-0.5,0.3) node {$+1$} ;
\draw (0.4,-0.4) node {\textcolor{oldgold}{$-1$}} ;
\draw (2.5,0.3) node {\textcolor{greyish}{$-1$}} ;
\draw (3.4,-0.4) node {\textcolor{oldgold}{$-1$}} ;
\draw (-0.5,3.3) node {$+1$} ;
\draw (0.4,2.6) node {\textcolor{blue}{$+1$}} ;
\draw (2.5,3.3) node {\textcolor{greyish}{$-1$}} ;
\draw (3.4,2.6) node {\textcolor{blue}{$+1$}} ;
\draw[fill=black,opaque] (0,-0.2) arc (270:90:0.2cm) -- cycle ;
\draw[fill=gold,opaque] (0,-0.2) arc (-90:90:0.2cm) -- cycle ;
\draw[fill=greyish,opaque] (3,-0.2) arc (270:90:0.2cm) -- cycle ;
\draw[fill=gold,opaque] (3,-0.2) arc (-90:90:0.2cm) -- cycle ;
\draw[fill=black,opaque] (0,2.8) arc (270:90:0.2cm) -- cycle ;
\draw[fill=blue,opaque] (0,2.8) arc (-90:90:0.2cm) -- cycle ;
\draw[fill=greyish,opaque] (3,2.8) arc (270:90:0.2cm) -- cycle ;
\draw[fill=blue,opaque] (3,2.8) arc (-90:90:0.2cm) -- cycle ;
\draw (1.5,-1.5) node  {\textit{Product} colours, $\sigma\tau$} ;
\end{tikzpicture}
\caption{Good (thick, black) and bad (thin, red) edges for the three Ising-type configurations
in the Ashkin-Teller model configuration in
Figure~\ref{fig:ashkin-teller-col}.}
%Figure~2.}
\label{fig:ashkin-teller-col-divns}
\end{center}
\end{figure}

These yield a third such function, their product $\sigma\tau\,:\,V\rightarrow\{\pm1\}$, given for
each $v\in V$ by $(\sigma\tau)(v)=\sigma(v)\tau(v)$.  So we have three Ising spin
functions, \textit{coupled} in the sense that they are not all independent: we can regard any two
of them as independent, but together they determine the third.  Each of them gives its own
division of $E(G)$ into good and bad edges, so we have three such divisions altogether:
\begin{center}
\begin{tabular}{ccccc}
spin   &  ~~~~  &   \multicolumn{3}{c}{edge type}   \\
  &&  \textbf{good}  & ~ & \textcolor{red}{\textbf{bad}}   \\   \hline
$\sigma$  &&  $E^-(\sigma)$  &&  $E^+(\sigma)$   \\
$\tau$  &&  $E^-(\tau)$  &&  $E^+(\tau)$   \\
$\sigma\tau$  &&  $E^-(\sigma\tau)$  &&  $E^+(\sigma\tau)$
\end{tabular}
\end{center}
We illustrate these divisions, for the graph and configuration of Figure~\ref{fig:ashkin-teller-col}, in Figure~\ref{fig:ashkin-teller-col-divns}.

The \textit{partition function of the \textit{symmetric} Ashkin-Teller model} is given by
\[
Z_{\hbox{\scriptsize SymAT}}(G;K,K') =
e^{(2K+K')|E|} ~~
\sum_{\sigma,\tau:V\rightarrow\{\pm1\}} ~~
e^{ - 2K \cdot |E^-(\sigma)|
-  2K \cdot |E^-(\tau)| 
-  2K' \cdot |E^-(\sigma\tau)| 
}.
\]
The symmetry is because the coefficients of $|E^-(\sigma)|$ and $|E^-(\tau)|$ in the exponent
are the same, namely $K$.
In the general four-state Ashkin-Teller model, these can be different, and the partition function depends
on three variables instead of the two ($K$, $K'$) we have here.

If $K'=0$ then $Z_{\hbox{\scriptsize SymAT}}(G;K,K')$ is just the square of an Ising model partition
function:
\[
Z_{\hbox{\scriptsize SymAT}}(G;K,0) = Z_{\scriptsize\hbox{Ising}}(G;K)^2.
\]
If $K=K'$ then $Z_{\hbox{\scriptsize SymAT}}(G;K,K')$ is just a Potts model partition function up
to a simple factor:
\[
Z_{\hbox{\scriptsize SymAT}}(G;K,K) = e^{-3K|E|} Z_{\scriptsize\hbox{Potts}}(G;4K,4).
\]

In these two cases, the symmetric Ashkin-Teller model partition function can be obtained from the
Tutte polynomial, since the Ising and Potts model partition functions can be so obtained.
But this is not possible in general.  Direct
evaluation shows that the two Tutte equivalent Gray graphs (Figure~\ref{fig:gray-graphs})
have different symmetric Ashkin-Teller model partition functions.
From this we see that the symmetric Ashkin-Teller model partition function is not a specialisation
of the Tutte polynomial.

\subsection{Interpolating between contraction and deletion}

We can represent any subset $X\subseteq E$ of the edge set of a graph by its
\textit{characteristic vector} $\mathbf{x}=(x_e)_{e\in E}\in\hbox{GF}(2)^E$ defined by
\[
x_e=\left\{
\begin{array}{cl}
1,  &  \hbox{if $e\in E$},   \\
0,  &  \hbox{if $e\not\in E$}.
\end{array}
\right.
\]
A \textit{cocircuit} of $G$ is a minimal set of edges whose removal does not disconnect
any component of $G$.  These are also the circuits of the dual of the cycle matroid of $G$.
The \textit{cocircuit space} of $G$ is the linear space over $\hbox{GF}(2)$ generated by
the characteristic vectors of cocircuits of $G$.  This is also the row space of the binary
incidence matrix of $G$, which has rows indexed by $V$, columns indexed by $E$, and
each entry is 1 or 0 according as its vertex is, or is not, incident with its edge.

In \cite{farr93b}, the operations of contraction and deletion are applied directly to the indicator
functions of cocircuit spaces.
Let $f : 2^E \rightarrow \{0,1\}$ be the indicator function of the cocircuit space of $G$.
Let $f\contract e: 2^E \rightarrow \{0,1\}$ and $f\delete e : 2^E \rightarrow \{0,1\}$ be the
indicator functions of $G/e$ and $G\setminus e$, respectively.  Then it is shown in \cite{farr93b} that
\[
(f\contract e)(X) = \frac{f(X)}{f(\emptyset)}\,, ~~~~~~
(f\delete e)(X) = \frac{f(X) + f(X\cup\{e\})}{
f(\emptyset) + f(\{e\})} \,.
\]
This is used in two related generalisations.
Firstly, in \cite{farr93b}, contraction and deletion are extended
to arbitrary $f : 2^E \rightarrow\mathbb{R}$ satisfying $f(\emptyset)=1$;
in later work it was convenient to include all
$f : 2^E \rightarrow\mathbb{C}$.  Such functions $f$ are called \textit{binary functions}, and if $f$
is the indicator function of the cocircuit space of a graph then it is \textit{graphic}.
Secondly, in \cite{farr04}, contraction and deletion are
considered to be just two specific reductions in a whole family
of $\lambda$\textit{-reductions}
parameterised by $\lambda\in\mathbb{R}$, and later by $\lambda\in\mathbb{C}$ in \cite{farr07b},
using the notation and expression in the middle column, where
$e\in E$ and $X\subseteq E\setminus\{e\}$:
\[
\begin{array}{ccccc}
\hbox{\textcolor{blue}{Contraction}} & ~~~~~ & 
\hbox{\textcolor{blue}{$\lambda$-reduction}} & ~~~~~ & 
\hbox{\textcolor{blue}{Deletion}}   \\
(\lambda=0)  &&&&  (\lambda=1)    \\
&&&&   \\
(f\contract e)(X)  &  & (f\minor{\lambda} e)(X)  &  &
(f\delete e)(X)   \\
&&&&   \\
{\displaystyle \frac{f(X)}{f(\emptyset)} }
&& \displaystyle \frac{f(X) + \lambda f(X\cup\{e\})}{
f(\emptyset) + \lambda f(\{e\})}
&& 
{\displaystyle \frac{f(X) + f(X\cup\{e\})}{f(\emptyset) + f(\{e\})} }
\end{array}
\]
Ordinary contraction and deletion are given by $\lambda=0$ and $\lambda=1$, respectively.
When $\lambda\in[0,1]$, one can think of the $\lambda$-reduction as interpolating
between contraction and deletion.  These $\lambda$-reductions come in dual pairs,
with the dual of $\lambda$-reduction being $\lambda^*$-reduction
where
\[
\lambda^*=\frac{1-\lambda}{1+\lambda} \,.
\]
When $\lambda\not\in\{0,1\}$, a $\lambda$-reduction of a graphic binary function is not,
in general, graphic.

These dual $\lambda$-reduction operations are accompanied by parameterised versions
of the rank function and Whitney rank generating function.  The $\lambda$\textit{-rank} function
$Q^{(\lambda)}f$ is defined by
\[
(Q^{(\lambda)}f)(X) = \log_2 \left(
\frac{(1 + \lambda^*)^{\left|V\right|} 
\sum_{W\subseteq E} \lambda^{\left|W\right|} f(W)}{
\sum_{W\subseteq E} \lambda^{\left|W\cap(E\setminus X)\right|} (
\lambda^*)^{\left|W\cap X\right|} f(W) }
\right) \,,
\]
where $X\subseteq E$.  The $\lambda$\textit{-Tutte-Whitney function} is defined by
\[
R_1^{(\lambda)}(G; x, y) = 
R_1^{(\lambda)}(f; x, y) =
y^{-Q^{(\lambda)}f(E)} \sum_{X\subseteq E}
(xy)^{Q^{(\lambda)} f(E) - Q^{(\lambda)} f(X)}
y^{\left|X\right|} .
\]
These functions are shown in \cite{farr04} to satisfy a generalisation of the contraction-deletion
relation of the same form as for Tutte-Whitney polynomials,
with $\lambda$-reductions and their dual $\lambda^*$-reductions taking the place
of contraction and deletion.  The \textit{loopiness} and \textit{coloopiness} of the element
$e\in S$ under the function $f$ are defined by the functions
\begin{eqnarray*}
\loopy^{(\lambda)}(f, e) & = &
Q^{(\lambda^*)}f(E) - Q^{(\lambda^*)}f(E\setminus e) ,   \\
\coloopy^{(\lambda)}(f, e) & = &
Q^{(\lambda)}f(E) - Q^{(\lambda)}f(E\setminus e)  ~~=~~  \loopy^{(\lambda^*)}(f, e) .
\end{eqnarray*}

\begin{thm}[\cite{farr04}]
\label{thm:lambda-del-con}
For any binary function $f:2^E\rightarrow\mathbb{C}$ and any $e\in E$,
\[
R^{(\lambda)}(f; x, y)  =  
x^{\scriptcoloopy^{(\lambda)}(f, e)} R^{(\lambda)}(f\minor{\lambda}e; x, y) +
y^{\scriptloopy^{(\lambda)}(f, e)} R^{(\lambda)}(f\minor{\lambda^*}e; x, y) .
\]
\end{thm}

\subsection{Go polynomials}
\label{sec:gopoly}

J\'anos is a keen Go player.  He and the first author have played several times;
as the stronger player, J\'anos plays with the white stones and gives a handicap.
So it is fitting to mention some surprising points
of contact between graph polynomials and the theory of Go.

Here we use the name by which the game is known in Japan and in the West.
It is called \textit{W\'eiq\'\i}\ in China and \textit{Baduk} in Korea.

Go is normally played on the square grid graph of $19\times19$ vertices.
%, see Figure \ref{fig:go19x19}.
But it has long been recognised that Go is an entirely graph-theoretic game and can be played
on any graph.  It has been played, for example,
on a three-dimensional diamond lattice graph \cite{segerman}
and on a map of Milton Keynes \cite{hunt2013}.

Thorpe and Walden \cite{thorp-walden64,thorp-walden72} seem to have been among the
first to formalise the rules and concepts of Go in order to support mathematical and
computational investigation.  Their work uses some graph-theoretic concepts but is still
embedded in rectangular grid graphs.  Benson \cite{benson76} uses graph-theoretic language
and concepts to characterise unconditional life of groups of stones.
Tromp and Taylor \cite{tromp-taylor,sensei} defined ``logical rules'' for Go, based on rules
from the New Zealand Go Society.  These rules are concise, precise, elegant, and purely
graph-theoretic.  Although they state initially that Go is played on a $19\times19$ square grid,
these rules may be used with only cosmetic modifications to play on any graph.

Let $G$ be a graph and $f$ be a partial $\lambda$-assignment of $G$.
A chromon of $(G,f)$ is \textit{free} if it has a vertex that is adjacent to an uncoloured vertex.
A partial $\lambda$-assignment $f$ is a \textit{legal $\lambda$-position} in $G$,
or just a \textit{legal position} if $\lambda$ is clear from the context,
if every chromon is free.

Normally, Go is played with just two colours, Black and White, so $\lambda=2$, although
multiplayer versions with $\lambda>2$ have been proposed and equipment (in the form
of coloured Go stones) is available for them.

In Go terminology, chromons correspond to what a Go player might call \textit{chains}.
A chromon is free if it is a chain with at least one liberty, and a legal position
is one in which every chain has a liberty.
When playing the game according to the rules, every position will be a legal position except
during the brief intervals after a capture is made and before all the captured stones are
removed from the board.  We do not discuss the rules of Go further; see, e.g.,~\cite{iwamoto72}
for more information.

Two \textit{Go polynomials} based on these concepts were introduced
in \cite{farr03}, studied further in \cite{farr-schmidt08,farr2017,farr2019b},
and mentioned by J\'anos in his first inventory of the Zoo \cite{makowsky06,makowsky08}.
One of them simply counts legal positions:
\[
\Go^{\#}(G;\lambda)  =
\hbox{number of legal $\lambda$-positions for $G$}.
\]
For example, it can be shown that
$\Go^{\#}(C_4; \lambda) = 1 + 14 \lambda^2$.
The other Go polynomial from \cite{farr03} is based on a simple probability model.
Let $p\le\frac{1}{2}$
and construct a random partial 2-assignment $f$ as follows.
Each $v\in V$ is assigned colour 1 or 2, with probability $p$ each, or is left
uncoloured with probability $r:=1-2p$, with decisions for different vertices being independent.
Under this model, define
\[
\Go(G;p)  =
\Pr(\hbox{$f$ is a legal 2-position for $G$}).
\]
We could also define a bivariate version of the second polynomial by extending
the probability model to $\lambda$ colours, where $f$ is now a partial $\lambda$-assignment,
the probability $p$ now satisfies $p\le\lambda^{-1}$, and $f$ is assigned colour
$k\in[\lambda]$ with probability $p$ for each colour
and is left uncoloured with probability $r:=1-\lambda p$:
\[
\Go(G;p,\lambda)  =
\Pr(\hbox{$f$ is a legal $\lambda$-position for $G$}).
\]
All three functions can be shown to be polynomials, so they can take as arguments any
$\lambda,p\in\mathbb{C}$ although we only know of combinatorial interpretations for the
values used in the definitions above.  This suggests the problem of finding combinatorial
interpretations at other values of $p$ and $\lambda$.  We suggest $\lambda=-1$ as one
that might be worth exploring, since the chromatic polynomial has an interesting combinatorial
interpretation at $\lambda=-1$, namely the number of acyclic orientations \cite{stanley73},
and Go polynomials can be expressed naturally as sums of chromatic polynomials \cite{farr03}.

These Go polynomials are all exponential-time computable, and they
seem unlikely to be polynomial-time computable because it was shown in \cite[\S4]{farr03}
that, for any fixed integer $\lambda\ge2$, computing the value of $\Go^{\#}(G;\lambda)$
for an input graph $G$  is \#P-hard, using methods from transcendental number theory.

\subsection{Polynomials for partial colourings}
\label{sec:polys-partial-cols}

Given a graph $G$, a nonnegative integer $\lambda$ representing some number of available
colours, and a probability $p\le\lambda^{-1}$, put $r:=1-\lambda p$ and define the following
random colouring model, noting that the partial $\lambda$-assignment it produces is not
necessarily a (proper) colouring.

Each vertex $v\in V(G)$ remains uncoloured with probability $r$ and otherwise is given a colour
chosen uniformly at random from the $\lambda$ available colours.  So, for any specific colour,
the probability that $v$ gets that colour is $p$.  The choices made at different vertices are
independent.  This process generates a random partial $\lambda$-assignment whose
domain is a subset of $V(G)$.

A partial $\lambda$-assignment $f$ of $G$ is a \textit{partial $\lambda$-colouring}
if it is a colouring of $G[\dom f]$, the subgraph of $G$ induced by the coloured vertices.
In a partial $\lambda$-colouring, vertices that are adjacent in $G$ cannot both get
the same colour; either they receive different colours or at least one of them is uncoloured.

We say $f$ is $\lambda$\textit{-extendable} if there is a $\lambda$-colouring of $G$ that extends $f$.

A vertex $v\not\in\dom f$ is \textit{immediately $\lambda$-forced by $f$} if its neighbours in $\dom f$
have exactly $\lambda-1$ distinct colours among them, and in this case we write $f;v$ for the
partial $\lambda$-assignment on $(\dom f)\cup\{v\}$ which agrees with $f$ on $\dom f$ and gives
$v$ the sole colour from $[\lambda]$
that does not appear among its neighbours.  The motivation is that, if vertex
$v$ is to receive a colour from $[\lambda]$ without creating any bad edges,
then it must be given that one colour that is unused by any of its neighbours.
If the number of different colours that appear among the neighbours of $v$ is $\lambda$, then
there is no hope for $v$: any colouring of it creates one or more bad edges.
If the neighbours of $v$ only have $\le\lambda-2$ colours,
then the colour of $v$ is not determined by the colours of its neighbours,
as there are at least two options for it.

If $f$ immediately $\lambda$-forces $v$, then there is only one possible colour for $v$ in all
$\lambda$-colourings of $G$ that extend $f$.  But the converse does not hold in general: the
colour of $v$ may be uniquely determined without being forced in this specific local sense.

The vertex $v$ is \textit{(eventually) $\lambda$-forced by $f$} if there is a sequence
of vertices $v_1,\ldots,v_k=v$, all outside $\dom f$,
and a sequence of partial $\lambda$-assignments $f=f_0,f_1,\ldots,f_k$ such that,
for all $i\in[\lambda]$,
\begin{itemize}
\item $\dom f_i=(\dom f)\cup\{v_1,\ldots,v_i\}$,
\item $v_i$ is immediately forced by $f_{i-1}$,
\item $f_i=f_{i-1};v_i$.
\end{itemize}
We may drop $\lambda$ from ``$\lambda$-forced'' when it is clear from the context.

Note that if $f$ is improper then it is not $\lambda$-extendable,
and that if $f$ is not $\lambda$-extendable then (regardless of whether or not it is proper)
it cannot force a $\lambda$-colouring of $G$.

For every graph $G$ we define three bivariate functions based on partial colourings.

The \textit{partial chromatic polynomial} $\hbox{PC}(G;p,\lambda)$ is defined by
\[
\hbox{PC}(G;p,\lambda)  =  \Pr(\hbox{$f$ is a $\lambda$-colouring of $G[\dom f]$}),
\]
where $p\in[0,\lambda^{-1}]$ is the probability that a specific vertex gets a particular colour,
and $\lambda\in\mathbb{N}$.
Since $\hbox{PC}(G;p,\lambda)$ is a polynomial (as we will shortly see), its domain is $\mathbb{C}^2$.
Observe that $\hbox{PC}(G;\lambda^{-1},\lambda)=\lambda^{-n}P(G;\lambda)$ and
$\hbox{PC}(G;0,\lambda)=1$.
The partial chromatic polynomial is a simple algebraic transformation of
the \textit{generalised chromatic polynomial} $P(G;x,y)$ introduced
by~\cite{dohmen-poenitz-tittmann03}.  For $x\in\mathbb{N}$ and $y\in[x]$, the value of $P(G;x,y)$
is the number of $[x]$-assignments for which the first $y$ colour classes are proper.
If we divide by the total number $x^n$ of all $[x]$-assignments, then this falls within
the random partial colouring framework by putting $p=x^{-1}$, $\lambda=y$ and $r=1-x^{-1}y$.
So we have
\[
P(G;x,y) ~=~ \frac{\hbox{PC}(G;x^{-1},y)}{x^n} \,.
\]
The change of variables is easily inverted to obtain    %  using  x = 1/p
\[
\hbox{PC}(G;p,\lambda) ~=~ \frac{P(G;p^{-1},\lambda)}{p^n} \,.
\]
By considering all possibilities for the domain of $f$, we obtain
\begin{equation}
\label{eq:partial-chrom-poly-as-sum}
\hbox{PC}(G;p,\lambda)  =  \sum_{C\subseteq V} P(G[C];\lambda)p^{|C|}(1-\lambda p)^{n-|C|},
\end{equation}
which is just a rearrangement of \cite[Theorem 1]{dohmen-poenitz-tittmann03} and shows that
$\hbox{PC}(G;p,\lambda)$ is indeed a polynomial.  When $\lambda$ is fixed, we denote the
resulting polynomial in $p$ by $\hbox{PC}_{\lambda}(G;p)$.

The partial chromatic polynomial is also a specialisation of some trivariate graph polynomials
that count edge-subsets according to $(|X|,\nu(X),\kappa(X))$
or according to simple invertible transformations of them, e.g., $(|X|,\nu(X),\rho(X))$.
Historically, the first of these trivariate polynomials seems to have been
the \textit{Borzacchini-Pulito polynomial} \cite{borzacchini-pulito82}, given by
\[
\hbox{BP}(G;x,y,z)
~=~ \sum_{X\subseteq E} x^{|X|} y^{k(X)} z^{\nu{X}}
~=~ \sum_{X\subseteq E} x^{|X|} y^{\kappa(X)+n-\nu(X)} z^{\nu{X}}
~=~ \sum_{X\subseteq E} x^{|X|} y^{n-\rho(X)} z^{\nu{X}} .
\]
They gave a ternary reduction relation for it in \cite[Theorem 2]{borzacchini-pulito82}.
Averbouch, Godlin and Makowsky \cite{averbouch-godlin-makowsky08,averbouch-godlin-makowsky2010}
introduced the \textit{edge elimination polynomial} $\xi(G;x,y,z)$ as the most general
trivariate graph polynomial that satisfies a ternary reduction relation using edge deletion, contraction
and extraction (deletion of an edge, its endpoints and their incident edges) and showed how
to express it as a sum over pairs of disjoint subsets of edges.  Trinks \cite{trinks2012} showed
that these two trivariate polynomials are equivalent in the sense that each can be transformed
to the other by simple transformations of the variables.

The \textit{extendable colouring function} is given by
\[
\hbox{EC}(G;p,\lambda)  =  \Pr(\hbox{$f$ is $\lambda$-extendable}) .
\]
Observe that $\hbox{EC}(G;\lambda^{-1},\lambda)=\lambda^{-n}P(G;\lambda)$ and
$\hbox{EC}(G;0,\lambda)$ is 1 if $G$ is 3-colourable and 0 otherwise.
We can express $\hbox{EC}(G;p,\lambda)$ as a sum over all possible domains for $f$,
\[
\hbox{EC}(G;p,\lambda)  =  \sum_{C\subseteq V} \hbox{EC}(G,C;\lambda) p^{|C|}(1-\lambda p)^{n-|C|},
\]
where $\hbox{EC}(G,C;\lambda)$ is the number of partial $\lambda$-colourings $f$ such that
$\dom f=C$ and $f$ has an extension which is a $\lambda$-colouring of $G$.
But this does not show that
$\hbox{EC}(G;p,\lambda)$ is a polynomial.  In fact, in general it is not.  Consider $K_2$.
\begin{equation}
\label{eq:EC-null}
\hbox{EC}(K_2;p,\lambda) = \left\{
\begin{array}{cl}
0,  &  \hbox{if $\lambda=1$},   \\
1-\lambda p^2,  &  \hbox{if $\lambda\ge2$}.
\end{array}
\right.
\end{equation}
It follows from (\ref{eq:EC-null}) that $\hbox{EC}(K_2;p,\lambda)$ is not a polynomial,
since when $p=0$ it is 1 for every positive integer $\lambda\ge2$ but is 0 for $\lambda=1$,
a property that no polynomial can have.  Nonetheless, $\hbox{EC}(G;p,\lambda)$ is a polynomial
in $p$ for every fixed $\lambda\in\mathbb{N}$.  When $\lambda$ is fixed, we denote this polynomial
in $p$ by $\hbox{EC}_{\lambda}(G;p)$.

For the extendable colouring polynomial $\hbox{EC}_{\lambda}(G;p)$, even checking the
structures being counted seems hard in general for $\lambda\ge3$.
Testing whether a given partial $\lambda$-assignment of a given graph $G$
is extendable to a $\lambda$-colouring of $G$ is NP-complete when $\lambda\ge3$,
by polynomial-time reduction from graph $\lambda$-colourability:
when the partial $\lambda$-assignment is the empty $\lambda$-assignment,
leaving all vertices uncoloured, we just have standard $\lambda$-colourability.
So this polynomial is likely to be more difficult to work with than many that have been studied.

The \textit{forced colouring function} is given by
\[
\hbox{FC}(G;p,\lambda)  =  \Pr(\hbox{$f$ eventually forces a $\lambda$-colouring of $G$}) .
\]
Again, we have $\hbox{FC}(G;\lambda^{-1},\lambda)=\lambda^{-n}P(G;\lambda)$.
We can express $\hbox{FC}(G;p,\lambda)$ as a sum,
\[
\hbox{FC}(G;p,\lambda)  =  \sum_{C\subseteq V} \hbox{FC}(G,C;\lambda) p^{|C|}(1-\lambda p)^{n-|C|},
\]
where $\hbox{FC}(G,C;\lambda)$ is the number of partial $\lambda$-colourings $f$ such that
$\dom f=C$ and $f$ eventually forces a $\lambda$-colouring of $G$.  But, again, we do not
have a polynomial, in general.  Consider a null graph.
\begin{equation}
\label{eq:FC-null}
\hbox{FC}(\overline{K_n};p,\lambda) = \left\{
\begin{array}{cl}
1,  &  \hbox{if $\lambda=1$},   \\
(\lambda p)^n,  &  \hbox{if $\lambda\ge2$}.
\end{array}
\right.
\end{equation}
This is because when there is only one colour,
the colour of every initially-uncoloured isolated vertex is forced,
while if there are at least two colours, an uncoloured isolated vertex cannot be forced, so it can
only get a colour from the initial random partial colouring.
It follows from (\ref{eq:FC-null}) that $\hbox{FC}(\overline{K_n};p,\lambda)$ is not a polynomial,
since when $p=0$ it is zero for every positive integer $\lambda\ge2$ but is 1 for $\lambda=1$.
But $\hbox{FC}(G;p,\lambda)$ is a polynomial in $p$
for every fixed $\lambda\in\mathbb{N}$, denoted by
$\hbox{FC}_{\lambda}(G;p)$.

Our three bivariate functions satisfy
\[
\hbox{FC}(G;p,\lambda) \le \hbox{EC}(G;p,\lambda) \le \hbox{PC}(G;p,\lambda)
\]
whenever $\lambda\in\mathbb{N}$ and $0\le p\le \lambda^{-1}$.

For these three functions, the case $\lambda=3$ is of particular interest:
\begin{eqnarray*}
\hbox{PC}_3(G;p)  & = &  \hbox{PC}(G;p,3),   \\
\hbox{EC}_3(G;p)  & = &  \hbox{EC}(G;p,3),   \\
\hbox{FC}_3(G;p)  & = &  \hbox{FC}(G;p,3).
\end{eqnarray*}
We are particularly interested in $\hbox{FC}_3(G;p)$ which we call the
\textit{forced 3-colouring polynomial} of $G$.  We list some basic examples and observations.
\begin{eqnarray*}
\hbox{FC}_3(\overline{K_n};p)  & = &  (3p)^n,   \\
\hbox{FC}_3(K_2;p)  & = &  6p^2,   \\
\hbox{FC}_3(K_3;p)  & = &  6p^2(3-8p),   \\
\hbox{FC}_3(K_{1,2};p)  & = &  6p^2(1-p),   \\
\hbox{FC}_3(G;{\textstyle\frac{1}{3}})  & = &  3^{-n} P(G;3),   \\
4^{n} \cdot\hbox{FC}_3(G;{\textstyle\frac{1}{4}})  & = &  \hbox{\# partial 3-assignments that eventually force a $\lambda$-colouring of $G$},   \\
\hbox{FC}_3(G\sqcup H;p)  & = &  \hbox{FC}_3(G;p)\hbox{FC}_3(H;p).
\end{eqnarray*}

In many graph polynomials, the structures being counted can be checked very quickly
\textit{in parallel}.  Typically, these checking problems belong to the class NC of decision
problems solvable by uniform families of logical circuits of polynomial size and polylogarithmic
depth, both of these being upper bounds in terms of the input size.
(See, e.g., \cite{greenlaw-hoover-ruzzo1995} for the theory of NC and P-completeness.)
For example, validity of a particular 3-colouring, or an independent set, or a clique, is
just a conjunction of local conditions based on the edges of the graph;
validity of a flow or a matching or a dominating set is a
conjunction of local conditions based on the vertex neighbourhoods.

One motivation for studying the forced 3-colouring polynomial is because of the computational
complexity of checking the structures being
counted, i.e., checking if a partial 3-assignment forces a 3-colouring of the graph.
This can be done in polynomial time, so the situation is not as bad as for extendable colouring
polynomials, and may seem more akin to classical graph polynomials based on enumerating
structures like colourings, independent sets, matchings and so on.
But the forced 3-colouring polynomial has the distinctive feature that its checking problem
is, in a precise sense, as hard as any problem in P.
Specifically, it is logspace-complete for P \cite{farr94}, which is considered to be evidence
that there is no fast parallel algorithm for this test and that it does not belong to NC.
Intuitively, this may make these graph polynomials more difficult to compute and to investigate
than the many others where the structure-checking can be done in NC.  So the study of them
might shed light on some less-explored regions of the theory of graph polynomials.

The forced 3-colouring polynomial does fall within the class of SOL-definable graph polynomials,
introduced by Makowsky and colleagues \cite{makowsky04,makowsky08} and explained in more
detail in \cite{godlin-katz-makowsky2012,kotek2012}.
To see this, we can express it as
\[
\hbox{FC}_3(G;p) ~=~ \sum_{(C_1,C_2,C_3)\in\mathcal{FC}_3(G)} p^{\sum_{i=1}^3 |C_i|}(1-3p)^{n-\sum_{i=1}^3 |C_i|},
\]
where $\mathcal{FC}_3(G)$ is the set of mutually disjoint triples $(C_1,C_2,C_3)$
whose corresponding partial 3-assignment $f$ ---
i.e., the partial 3-assignment $f$ defined by $C_f(i)=C_i$ for $i\in\{1,2,3\}$ ---
forces a 3-colouring of $G$.  The condition that $(C_1,C_2,C_3)\in\mathcal{FC}_3(G)$
can be expressed in SOL using the following observation.

\begin{propn}
A partial $\lambda$-assignment $f$ forces a $\lambda$-colouring of $G$ if and only if
it has no extension $g$ which is not total and forces no vertex outside $\dom g$.
\end{propn}

\pf
($\Longrightarrow$)
Suppose $f$ forces a $\lambda$-colouring of $g$.  Then every extension of $f$ forces
the same $\lambda$-colouring of $G$.  So there is no non-total extension of $f$ that does not force
any vertex.

($\Longleftarrow$)
If $f$ does not force a $\lambda$-colouring of $G$, then the forcing process must stop with
at least one vertex uncoloured.  When that happens, the partial $\lambda$-assignment that
has been forced so far is an extension of $f$ that is not total and forces no vertex outside its domain.
\eopf

This observation justifies the following logical formulation of $\mathcal{FC}_3(G)$, which is
similar in design to the SOL expression for connectedness in \cite[\S3.3]{godlin-katz-makowsky2012}.
\begin{eqnarray*}
\lefteqn{(C_1,C_2,C_3)\in\mathcal{FC}_3(G)
~\Longleftrightarrow~}   \\
\lefteqn{\neg\exists (D_1,D_2,D_3) :}  \\
&&
\left( \bigwedge_{i=1}^3 (C_i\subseteq D_i)\right)  \wedge
\left(\bigwedge_{i=1}^3 \bigwedge_{j=i+1}^3 (D_i\cap D_j = \emptyset )\right) \wedge
\left( D_1\cup D_2\cup D_3\not=E \right)   \wedge   \\
&&
\forall v\in E\setminus(D_1\cup D_2\cup D_3) :   \\
&&   ~~
\left( \,
\left( \bigwedge_{i=1}^3 (\exists w\in D_i : vw\in E ) \right)
\vee
\bigvee_{i=1}^3 \bigvee_{j=i+1}^3 \left(\forall w\in V :  vw\in E \rightarrow \left(w\not\in D_i\cup D_j\right)\right)
\right) .
\end{eqnarray*}
This can all be expressed using the formalism of SOL,
with sets of vertices represented as unary relations,
adjacency as a binary relation, and so on.

We study the forced colouring function of a graph further, along with its associated polynomials, 
in \cite{farr2024}.

\section{Origins of graph polynomials}
\label{sec:origins}

Having discussed a variety of graph polynomials in the previous section,
it is a good time to review the origins of graph polynomials in general.

Historically, graph polynomials have been created in several different ways.
\begin{enumerate}[(i)]
\item
Sometimes, sequences of numbers that count structures of interest are used as coefficients
to define a generating polynomial for them.
For example, polynomials that count
independent sets \cite{helgason74,gutman-harary83},
cliques \cite{farrell89,zykov1949,zykov1964},
matchings \cite{heilmann-lieb70,farrell79,godsil-gutman81} and
dominating sets \cite{arocha-llano00} arose in this way.
\item
Sometimes (but less commonly), a sequence of numbers that count structures of interest
in a graph is taken to give the values of a function of the graph and an integer parameter,
and the function turns out to be a polynomial in that parameter.
This is how the chromatic polynomial arose~\cite{birkhoff1912-1913}.
\item
Sometimes, a probability model is defined on graphs, with independent, identically-distributed
random choices being made throughout the graph (typically on vertices or edges) according to
some probability $p$, and the probability that this gives a particular type of structure is a polynomial
in $p$.  Polynomials defined this way include
the all-terminal network reliability polynomial \cite{vanslyke-frank1971} and the stability
polynomial \cite{helgason74,farr93a}, as well as the polynomials based on partial colourings
that we introduced in \S\ref{sec:polys-partial-cols}.
Such polynomials are often easily transformed into generating
polynomials, as is the case for the relationship between the stability and independent set
polynomials.
\item
Sometimes, a graph invariant of physical interest is defined which turns out to be a polynomial,
possibly after an algebraic change of variables, as in the case of the partition functions of
the Ising model \cite{ising1925}, Ashkin-Teller model \cite{ashkin-teller43} and
Potts model \cite{potts52}.  There is some overlap between this type and some previous types,
e.g., the matching polynomial is framed as a partition function in \cite{heilmann-lieb70},
and the reliability polynomial models the survivability of networks in the presence of local link
failures.
\item
\label{item:tutte-poly}
Sometimes, certain graph invariants
are found to satisfy reduction relations of some kind,
and this motivates the development of a common generalisation which turns out to
be a polynomial.  This is how the Tutte polynomial was created,
as Tutte relates in \cite[p.~53]{tutte98} and \cite{tutte04}.
\item
Sometimes, a polynomial is created by analogy with existing polynomials and/or by
generalising them.  See \cite[\S3.4]{farr07a} for a discussion of the different ways in which
Tutte-Whitney polynomials have been generalised; the informal classification given there could
apply to analogues and generalisations of any graph polynomial.
An early example of this is the Whitney rank generating function itself, which arose in 1932 as
a bivariate generalisation of the chromatic polynomial \cite{whitney32d} without noting any other
combinatorial interpretations of its values and with the deletion-contraction relation only being
a ``note added in proof'' and attributed to R. M. Foster.
The Oxley-Whittle polynomial \cite{oxley-whittle93a,oxley-whittle93b} arose as the analogue
of Tutte-Whitney polynomials for 2-polymatroids.  Tutte-Whitney polynomials of graphs were also the
inspiration for the various topological Tutte polynomials for ribbon graphs and embedded graphs
\cite{chmutov2022,ellis-monaghan-moffatt2013}.
\item
Occasionally, a polynomial is defined by specialisation from an existing polynomial, instead of
by generalisation: it is not always the case that the particular motivates the general.
The flow polynomial was first discovered (but not named) by Whitney \cite{whitney33}
as the dual of the chromatic polynomial and a specialisation of his rank generating function,
but without giving its values any combinatorial interpretation; see \cite[\S34.7]{farr2022}.
Only later did Tutte identify the connection with flows \cite{tutte54}.
\item
When counting graphs and other combinatorial objects up to symmetry,
cycle index polynomials are used: see \cite{harary-palmer73}
for their use in counting various types of graphs up to isomorphism, a field which originated
with Redfield \cite{redfield1927} and P\'olya \cite{polya1937}.
Beginning with work by Cameron and collaborators
in the early 2000s, cycle index polynomials have been used to define graph polynomials
that count colourings and other structures up to symmetries under the automorphism group
of the graph \cite{cameron02,cameron07,cameron-jackson-rudd08}.
\item
Sometimes, the characteristic polynomial of a square matrix associated with a graph is studied,
the main example being
the characteristic polynomial of a graph\footnote{not to be confused with the characteristic polynomial of a matroid, which generalises the chromatic polynomial of a graph.}
which is obtained from the adjacency matrix.
\item
Sometimes, a polynomial is defined for other mathematical objects, and a natural construction
of graphs from those objects leads to a graph polynomial.   Some knot
polynomials have been translated to graph polynomials and linked to the Tutte polynomial:
see, e.g., \cite{huggett2022}.  The weight enumerator of a linear code over GF(2)
may be regarded as a binary matroid polynomial whose coefficients count members of a
circuit space according to their size (and similarly for a cocircuit space), thereby yielding a graph
polynomial through the special case of graphic matroids (see, e.g., \cite[\S15.7]{welsh1976}).
\end{enumerate}
This classification is just a set of observations of how the field has developed so far and
is certainly not any kind of prescription for how it must develop in the future.
It may not be exhaustive, and the categories may overlap.  Indeed it is common for a graph
polynomial to have multiple different formulations, so that it can be \textit{defined} in two or more
of the above ways, even if it was historically \textit{created} in just one way.
There is no single best route to defining graph polynomials: all the routes listed above have led
to new polynomials that have generated a lot of interest and some profound research.

Nonetheless, it is worth reflecting on the various inspirations for the diverse graph polynomials
that have been developed.
Analogy and generalisation have been very fruitful, but it is arguable that the real worth of a new graph
polynomial lies not so much in how well its theory echoes those of other known polynomials, but rather
in the information it contains about the graph (including the relationships it reveals between different features of it) and in the accessibility of that information.

In this context, the originality of Tutte's approach (\ref{item:tutte-poly}) is striking.
It did not really fall within any earlier approach.  It was grounded in important graph invariants.
His polynomial emerged naturally as a framework that captured what those invariants had
in common.  The fact that it is a \textit{polynomial} was a happy outcome, and not surprising since
one of the invariants he was abstracting from was the chromatic polynomial, but it does not seem
to have been an objective in itself.

It is conceivable that other collections of enumerative (or probabilistic) graph invariants
are waiting to be brought into common frameworks, and that these frameworks may sometimes
be quite different to Tutte's, and may not always lead to polynomials.

\section{Reduction relations}
\label{sec:redn-relns}

The recursive relation (\ref{eq:tgxy}) for the Tutte polynomial has the following characteristics.
\begin{itemize}
\item The number of cases is fixed (i.e., independent of the size of the graph).
\item In the base case, the expression is a fixed polynomial (in this case, just the constant 1).
\item In each other case,
the expression is linear (treating the Tutte polynomial
symbols $T(G;x,y)$, $T(G\setminus e;x,y)$, $T(G/e;x,y)$ as indeterminates,
with coefficients being fixed polynomials in $x,y$).
\item The graphs appearing in the right-hand sides are obtained from
the original graph $G$ by simple local operations.
\item Certain types of edges must be treated as special cases (namely, loops and coloops).
Such edges are ``degenerate'' in some sense, and the expression on the right-hand side is usually
simpler than in the general case, with fewer terms.
\item The number of terms in each of these linear expressions is fixed.
\item The polynomial is well defined, in that the order in which local operations are used
when evaluating the polynomial using (\ref{eq:tgxy}) does not matter \cite{tutte54}.
\end{itemize}

Many other graph polynomials also satisfy \textit{reduction relations} of
this type.  Examples include the independent set, stability, clique
and matching polynomials, the Oxley-Whittle polynomial for graphic
2-polymatroids \cite{oxley-whittle93a,oxley-whittle93b}, the Borzacchini-Pulito polynomial \cite{borzacchini-pulito82}, and the rich class of Tutte polynomials of Hopf algebras \cite{krajewski-moffatt-tanasa2018}.\footnote{Some notable
graph polynomials are simple transformations or partial evaluations of the Tutte polynomial,
so they satisfy reduction relations of this type because the Tutte polynomial does.
These include
the chromatic, flow, reliability, Martin and Jones polynomials,
the Ising and Potts model partition functions, the Whitney rank generating function
and the coboundary polynomial.}
One of the earliest studies of graph polynomials defined by a variety of reduction relations
beyond deletion-contraction is due to Zykov \cite{zykov1964};
this has roots in his earliest work in graph theory \cite{zykov1949}.

Godlin, Katz and Makowsky have given a general definition of reduction relations of roughly
the above type,
in the context of the logical theory of graph polynomials \cite{godlin-katz-makowsky2012}.
Paraphrasing, they showed that every graph polynomial with such a recursive definition in SOL can
be expressed as a SOL-definable sum over subsets.  This links, in one direction, two of
the main ways of defining graph polynomials: reduction relations, and sums over subsets (``subset
expansions" in their terminology).   They ask whether the link goes the other way too:
does every graph polynomial expressible as a SOL-definable sum over subsets satisfy a
reduction relation of their general SOL-based form?

Some graph polynomials do not seem to have a natural reduction relation within the class of
objects over which they are initially defined.
But it often happens that, even in such cases, there is a wider class of objects
to which the graph polynomial can be generalised and which
supports a reduction relation of this type.  We illustrate this with some examples.

\subsection{Counting edge-colourings}
\label{sec:redn-relns-edge-col}

For any graph $G$, define $P'(G;q)$ to be the number of $q$-edge-colourings
of $G$ (mirroring the standard use of $\chi'(G)$ and $\chi(G)$ to denote the chromatic index\footnote{i.e., the minimum $q$ such that $G$ is $q$-edge-colourable}
and chromatic number, respectively).
It is well known that edge-colourings of a graph $G$ correspond to
vertex-colourings of the line graph $L(G)$ of $G$ (see, e.g., \cite{fiorini-wilson77}).
So we have
\begin{equation}
\label{eq:edge-col-poly-chrom-poly-line-gph}
P'(G;q) = P(L(G);q),
\end{equation}
which makes plain that $P'(G;q)$ a polynomial.

The most natural way to get a reduction relation for $P'(G;q)$
is to use (\ref{eq:edge-col-poly-chrom-poly-line-gph}) and invoke the deletion-contraction
relation for chromatic polynomials:\footnote{We assume that line graphs have no loops, which
is the usual practice.  If, instead, we assume that loops are created in the line graph
corresponding to loops in $G$, then it is natural to assume that a line graph with a loop has no
edge-colouring.  If we were to allow colouring of loops in edge-colourings,
then our reduction relation would not work.}
\[
P(L(G);q)  ~=~ \left\{
\begin{array}{ll}
q^m,  &  \hbox{if $L(G)$ has no edges;}   \\
P(L(G)\setminus e;q) - P(L(G)/e;q),~~~~  &  \hbox{otherwise.}   \\
\end{array}
\right.
\]
But this relation holds in the class of \textit{all} graphs since, in general,
$L(G)\setminus e$ and $L(G)/e$ are not line graphs.

The base cases for this reduction relation are graphs with no edges.  Such a graph is a line graph
of a disjoint union of a matching and a set of isolated vertices.

\subsection{The symmetric Ashkin-Teller model}
\label{sec:redn-relns-symAT}

In \S\ref{sec:partn-fns}
we saw that the symmetric Ashkin-Teller model partition function is not a specialisation
of the Tutte polynomial.  It therefore does not obey the
kind of deletion-contraction relation characteristic of evaluations of the Tutte polynomial.

However, it turns out that it does satisfy such a relation in the wider class of binary functions.
This follows from the following result in \cite{farr07b} which shows that it is a specialisation
of a suitable $\lambda$-Tutte-Whitney function.

\begin{thm}[\cite{farr07b}]
\label{thm:ashkin-teller-lambda-tutte-whitney}
For each $K$ and $K_{\sigma\tau}$ there exists $\lambda$
such that the partition function $Z(G;K,K_{\sigma\tau})$
of the symmetric Ashkin-Teller model on a graph $G$
can be obtained from the $\lambda$-Tutte-Whitney function.
\end{thm}

Details, including expressions for $\lambda$ in terms of $K$ and $K_{\sigma\tau}$, are given in \cite{farr07b}.

Putting Theorems~\ref{thm:lambda-del-con} and~\ref{thm:ashkin-teller-lambda-tutte-whitney}
together gives a reduction relation for the Ashkin-Teller
model partition function in the class of binary functions.

\subsection{Go polynomials}
\label{sec:redn-reln-gopoly}

None of the Go polynomials introduced in \S\ref{sec:gopoly} and \cite{farr03} have an obvious
recurrence relation, of the type we have been considering, on graphs.
But there are recurrence relations in a wider class of objects that generalise graphs.
In \cite[\S3]{farr03}, it is found that $\Go^{\#}(G;\lambda)$ and $\Go(G;p)$ satisfy a family of
linear recurrence relations over a larger class of graphs, there called $\mathcal{L}$-graphs,
in which graphs may have extra labels on some of their vertices and edges that modify the
conditions that a partial $\lambda$-assignment must satisfy in order to be a legal position.

\subsection{Polynomials based on partial colourings}
\label{sec:redn-reln-partial-col}

We now consider reduction relations for the polynomials we introduced in \S\ref{sec:polys-partial-cols}.

The partial chromatic polynomial does not obey a deletion-contraction relation of the usual type,
as it is not obtainable from the Tutte polynomial and in fact contains extra information.
No reduction relation for it is given explicitly in \cite{dohmen-poenitz-tittmann03}.
But (\ref{eq:partial-chrom-poly-as-sum}) points to a reduction relation based on labelling.
A \textit{chromatically labelled graph} is a graph in which
each vertex may be labelled
\texttt{C}, indicating that it must receive a \underline{c}olour,
or \texttt{U}, indicating that it must be \underline{u}ncoloured; a vertex may have no label, but it
cannot have two labels.  Each chromatically labelled graph is written as $G^{(C,U)}$ where $C, U\subseteq V$ and $C\cap U=\emptyset$.
A \textit{totally chromatically labelled graph} $G^{(C,U)}$ is a chromatically labelled graph in which
every vertex is labelled, i.e., $C\cup U=V$.

A \textit{partial $\lambda$-assignment} of a chromatically labelled graph $G^{(C,U)}$
is a partial $\lambda$-assignment $f$ of $G$ such that $C\subseteq\dom f\subseteq V\setminus U$.
It is a \textit{partial $\lambda$-colouring} of $G^{(C,U)}$ if it is also a partial $\lambda$-colouring of $G$.
 
For chromatically labelled graphs, put
\begin{eqnarray*}
\hbox{PC}(G^{(C,U)};p,\lambda)   & = &  \Pr(\hbox{$f$ is a partial $\lambda$-colouring of $G^{(C,U)}$})   \\
& = &  \Pr(\left(\hbox{$f$ is a $\lambda$-colouring of $G[\dom f]$})\wedge(C\subseteq\dom f\subseteq V\setminus U)\right).
\end{eqnarray*}
This polynomial, in this wider class, has a simple reduction relation.

\begin{thm}
\label{thm:redn-reln-partial-chrom-poly}
For any $v\in V\setminus (C\cup U)$,
\[
\hbox{\rm PC}(G^{(C,U)};p,\lambda) ~=~
\hbox{\rm PC}(G^{(C\cup\{v\},U)};p,\lambda) +
\hbox{\rm PC}(G^{(C,U\cup\{v\})};p,\lambda) .
\]
\end{thm}

\pf
\begin{eqnarray*}
\lefteqn{\hbox{PC}(G^{(C,U)};p,\lambda)}   \\
& = &
\Pr((\hbox{$f$ is a $\lambda$-colouring of $G[\dom f]$})\wedge
(C\subseteq\dom f\subseteq V\setminus U))   \\
& = &
\Pr((\hbox{$f$ is a $\lambda$-colouring of $G[\dom f]$})\wedge
(C\subseteq\dom f\subseteq V\setminus U)\wedge(v\in\dom f))     \\
&&   +
\Pr((\hbox{$f$ is a $\lambda$-colouring of $G[\dom f]$})\wedge
(C\subseteq\dom f\subseteq V\setminus U)\wedge(v\not\in\dom f))   \\
& = &
\Pr((\hbox{$f$ is a $\lambda$-colouring of $G[\dom f]$})\wedge
(C\cup\{v\}\subseteq\dom f\subseteq V\setminus U))     \\
&&   +
\Pr((\hbox{$f$ is a $\lambda$-colouring of $G[\dom f]$})\wedge
(C\subseteq\dom f\subseteq V\setminus (U\cup\{v\})))   \\
& = &
\hbox{PC}(G^{(C\cup\{v\},U)};p,\lambda) +
\hbox{PC}(G^{(C,U\cup\{v\})};p,\lambda) .
\end{eqnarray*}
\eopf

The reduction relation of Theorem \ref{thm:redn-reln-partial-chrom-poly} cannot be used
on totally chromatically labelled graphs, when $C\cup U=V$.  In that case,
the partial chromatic polynomial is just a scaled version of
chromatic polynomial of $G-U$:
\begin{equation}
\label{eq:partial-chrom-poly-totally-chrom-labelled-gph}
\hbox{PC}(G^{(V\setminus U,U)};p,\lambda) = (1-\lambda p)^{|U|}p^{n-|U|} \, P(G-U;\lambda) \,.
\end{equation}
We return to this point in \S\ref{sec:levels}.

It is also possible to get a reduction relation for $\hbox{PC}(G;p,\lambda)$ on certain
vertex-weighted graphs using the fact that the partial chromatic polynomial is a specialisation of the $U$-polynomial of
Noble and Welsh \cite{noble-welsh1999} (because its equivalent polynomial $\xi(G;x,y,z)$ is)
which has a reduction relation on those weighted graphs.

We now consider extendable colouring polynomials and
show that $\hbox{EC}_{\lambda}(G)$ satisfies a reduction relation
in the class of chromatically labelled graphs we introduced above.
A partial $\lambda$-assignment of $G^{(C,U)}$ is $\lambda$-\textit{extendable} in $G^{(C,U)}$
if, as a partial $\lambda$-assignment of $G$, it is $\lambda$-extendable.
So, although a label \texttt{U} on a vertex specifies that it is uncoloured by our random
$\lambda$-assignment $f$, the vertex is allowed to be coloured by an extension of $f$.

For chromatically labelled graphs, put
\begin{eqnarray*}
\hbox{EC}(G^{(C,U)};p,\lambda)   & = &  \Pr(\hbox{$f$ is $\lambda$-extendable in $G^{(C,U)}$})   \\
& = &  \Pr\left((\hbox{$f$ is $\lambda$-extendable in $G$})\wedge(C\subseteq\dom f\subseteq V\setminus U)\right).
\end{eqnarray*}

\begin{thm}
\label{thm:redn-reln-ext-col-poly}
For any $v\in V\setminus (C\cup U)$,
\[
\hbox{\rm EC}_{\lambda}(G^{(C,U)};p)   ~=~
\lambda p \, \hbox{\rm EC}_{\lambda}(G^{(C\cup\{v\},U)};p) +
(1-\lambda p) \, \hbox{\rm EC}_{\lambda}(G^{(C,U\cup\{v\})};p) .
\]
\end{thm}

\pf
\begin{eqnarray*}
\hbox{EC}_{\lambda}(G^{(C,U)};p)
& = &
\Pr(\hbox{$f$ is $\lambda$-extendable in $G^{(C,U)}$})   \\
& = &
\Pr(\hbox{$f$ is $\lambda$-extendable in $G^{(C,U)}$}\mid v\in\dom f) \Pr(v\in\dom f) +   \\
&&
\Pr(\hbox{$f$ is $\lambda$-extendable in $G^{(C,U)}$}\mid v\not\in\dom f) \Pr(v\not\in\dom f)   \\
& = &
\lambda p \Pr(\hbox{$f$ is $\lambda$-extendable in $G^{(C\cup\{v\},U)}$}) +   \\
&&   (1-\lambda p) \Pr(\hbox{$f$ is $\lambda$-extendable in $G^{(C,U\cup\{v\})}$})   \\
& = &
\lambda p \, \hbox{EC}_{\lambda}(G^{(C\cup\{v\},U)};p) +
(1-\lambda p) \, \hbox{EC}_{\lambda}(G^{(C,U\cup\{v\})};p) .
\end{eqnarray*}
\eopf

Lastly, we consider forced colouring polynomials.
A partial $\lambda$-assignment of $G^{(C,U)}$ \textit{forces} a $\lambda$-colouring of $G^{(C,U)}$
if, as a partial $\lambda$-assignment of $G$, it forces a $\lambda$-colouring of $G$.
Again, a label \texttt{U} on a vertex only prevents it from being coloured by $f$ itself;
the label does not stop it from being eventually forced by $f$.
Put
\begin{eqnarray*}
\hbox{FC}(G^{(C,U)};p,\lambda)   & = &  \Pr(\hbox{$f$ forces a $\lambda$-colouring of $G^{(C,U)}$})   \\
& = &  \Pr\left((\hbox{$f$ forces a $\lambda$-colouring of $G$})\wedge(C\subseteq\dom f\subseteq V\setminus U)\right).
\end{eqnarray*}

A very similar argument
to Theorem \ref{thm:redn-reln-ext-col-poly} shows that
$\hbox{FC}_{\lambda}(G;p)$ satisfies a reduction relation
in the class of chromatically labelled graphs.

\begin{thm}
\label{thm:redn-reln-forced-col-poly}
For any $v\in V\setminus (C\cup U)$,
\[
\hbox{\rm FC}_{\lambda}(G^{(C,U)};p)   ~=~
\lambda p \, \hbox{\rm FC}_{\lambda}(G^{(C\cup\{v\},U)};p) +
(1-\lambda p) \, \hbox{\rm FC}_{\lambda}(G^{(C,U\cup\{v\})};p) .
\]
\eopf
\end{thm}

%\pf
%\begin{eqnarray*}
%\hbox{FC}_{\lambda}(G^{(C,U)};p)
%& = &
%\Pr(\hbox{$f$ forces a $\lambda$-colouring of $G^{(C,U)}$})   \\
%& = &
%\Pr(\hbox{$f$ forces a $\lambda$-colouring of $G^{(C,U)}$}\mid v\in\dom f) \Pr(v\in\dom f) +   \\
%&&
%\Pr(\hbox{$f$ forces a $\lambda$-colouring of $G^{(C,U)}$}\mid v\not\in\dom f) \Pr(v\not\in\dom f)   \\
%& = &
%\lambda p \Pr(\hbox{$f$ forces a $\lambda$-colouring of $G^{(C\cup\{v\},U)}$}) +   \\
%&&   (1-\lambda p) \Pr(\hbox{$f$ forces a $\lambda$-colouring of $G^{(C,U\cup\{v\})}$})   \\
%& = &
%\lambda p \, \hbox{FC}_{\lambda}(G^{(C\cup\{v\},U)};p) +
%(1-\lambda p) \, \hbox{FC}_{\lambda}(G^{(C,U\cup\{v\})};p) .
%\end{eqnarray*}
%\eopf

\subsection{Questions}

We have now seen six examples of graph polynomials which do not seem to satisfy a local linear
reduction relation over the class of graphs
but which do satisfy such relations over some larger class.

Other examples exist.
The $U$-polynomial of a graph, introduced by Noble and Welsh,
has a reduction relation in the larger class of vertex-weighted graphs in which the weights are positive
integers \cite{noble-welsh1999} (see also \cite{noble2022}).
Krajewski, Moffatt and Tanasa \cite{krajewski-moffatt-tanasa2018}
used their Hopf algebra framework to show that various topological Tutte polynomials without
full reduction relations\footnote{because the known relations did not cover all possible edge types}
can be extended to a larger class of objects (by augmenting the embedded graphs with some extra structure) so that they do have full reduction relations.  (See especially \cite[Remark 62]{krajewski-moffatt-tanasa2018}.)

These examples raise the question of how widespread this phenomenon is.
Which graph polynomials exhibit this phenomenon?  Can they be characterised in some
formal, rigorous way?
For the polynomials considered in \S\ref{sec:redn-reln-gopoly} and \S\ref{sec:redn-reln-partial-col},
suitable superclasses of the class of graphs can be found by introducing new labels on vertices
and/or edges with specific technical meanings.  This is likely to be a wider phenomenon and
may be able to be captured using the logical framework of Makowsky and colleagues.
But our first two cases, in \S\ref{sec:redn-relns-edge-col} and \S\ref{sec:redn-relns-symAT},
are not of this type, and it is not clear how to include them in a general characterisation of this
phenomenon.

\section{Levels of recursion}
\label{sec:levels}

For many graph polynomials, repeated application of a reduction relation leaves only trivial
graphs.  For example, repeated application of deletion-contraction relations for the chromatic
or Tutte polynomials leaves null graphs.  But sometimes a graph
polynomial has a reduction relation in which the base cases are themselves nontrivial graphs
and another reduction relation needs to be applied in order to reduce them to simpler
base cases; we might say that we have two ``levels'' of reduction relation.

Partial chromatic polynomials provide an example.  The reduction relation we gave
in Theorem~\ref{thm:redn-reln-partial-chrom-poly}
may be used to reduce the partial chromatic polynomial to a sum involving
partial chromatic polynomials of totally chromatically labelled graphs,
and each of these polynomials can in turn be expressed in terms of a chromatic polynomial,
by~(\ref{eq:partial-chrom-poly-totally-chrom-labelled-gph}).
So we can apply the deletion-contraction relation to each of the chromatic polynomials,
thereby expressing the partial chromatic polynomial as a sum of simple base cases --- chromatic
polynomials of null graphs --- with two levels of reduction.

It can get worse than this!
Go polynomials (\S\ref{sec:gopoly}) may be regarded as having three levels of
reduction.
Firstly, graphs are reduced to $\mathcal{L}$-graphs using \cite[Cor.~6 or Theorem 7]{farr03}.
Then $\mathcal{L}$-graphs are reduced to ordinary graphs again using \cite[Cor.~10]{farr03}.
Finally, these ordinary graphs are reduced to null graphs using the deletion-contraction relation
for the chromatic polynomial.  The paper in fact gives a method of expressing a Go polynomial
as a large sum of chromatic polynomials.

For some polynomials, the situation is less clear.  For extendable colouring and forced
colouring polynomials, we gave reduction relations in Theorems \ref{thm:redn-reln-ext-col-poly}
and \ref{thm:redn-reln-forced-col-poly} whose base cases require computation of those
polynomials for totally chromatically labelled graphs.  Those computations do not have
obvious analogues of~(\ref{eq:partial-chrom-poly-totally-chrom-labelled-gph}).

For extendable colouring polynomials, here is an attempt involving an addition-identification
relation on the set of vertices labelled \texttt{C}:
\[
\hbox{EC}_{\lambda}(G^{(C,U)};p) ~=~ 
\hbox{EC}_{\lambda}(G^{(C,U)}+uv;p) +
\hbox{EC}_{\lambda}(G^{(C,U)}/uv;p),
\]
for any $u,v\in C$ such that $uv\not\in E(G)$.  This can be applied repeatedly until the vertices
labelled \texttt{C} form a clique in $G$.  Whenever this clique has $>\lambda$ vertices, the
polynomial is identically 0 since the graph is not $\lambda$-colourable.  So we end up with
a sum of polynomials $\hbox{EC}_{\lambda}(H^{(D,U)};p)$ of totally chromatically labelled
graphs of the form $H^{(D,U)}$ where $D$ is a clique of size $\le\lambda$ in $H$.
Because $D$ is a clique, a partial $\lambda$-assignment of $H^{(D,U)}$ is $\lambda$-extendable
if and only if $H$ is $\lambda$-colourable.  So we have
\begin{eqnarray*}
\lefteqn{\hbox{EC}_{\lambda}(H^{(D,U)};p)}   \\
& = &
\Pr\left( (\hbox{$f$ is $\lambda$-extendable in $H$})\wedge(D\subseteq\dom f\subseteq 
V\setminus U)\wedge(\hbox{$f$ is a $\lambda$-colouring of $H[D]$}) \right)   \\
& = &
\Pr\left( (D\subseteq\dom f\subseteq V\setminus U)\wedge(\hbox{$f$ is a $\lambda$-colouring of $H[D]$}) \right) \times   \\
&&
\Pr\left((\hbox{$f$ is $\lambda$-extendable in $H$})\mid(D\subseteq\dom f\subseteq 
V\setminus U)\wedge(\hbox{$f$ is a $\lambda$-colouring of $H[D]$})\right)   \\
& = &
%\lambda p (\lambda-1) p \cdots p
p^{|D|}(\lambda)_{|D|} 
\cdot
\llbracket\hbox{$H$ is $\lambda$-colourable}\,\rrbracket.
\end{eqnarray*}
Here we have used the Iverson bracket:
\[
\llbracket\hbox{$H$ is $\lambda$-colourable}\,\rrbracket ~=~
\left\{
\begin{array}{cl}
1,  &  \hbox{if $H$ is $\lambda$-colourable};   \\
0,  &  \hbox{otherwise}.
\end{array}
\right.
\]
We seem to have gained something, computationally, by expressing it this way: we ``only'' have an
NP-complete quantity to evaluate, rather than a \#P-hard graph polynomial!  But it is no longer
a natural sum of graph polynomials, so it does not give us another layer of reduction relations
of the kind we have been considering.

We can try a similar approach with forced colouring polynomials
(Theorem~\ref{thm:redn-reln-forced-col-poly}).
Again, the base cases for the first reduction relation we use
are totally chromatically labelled graphs,
and again we can use addition-identification repeatedly to get a sum over totally
chromatically labelled graphs $H^{(D,U)}$
in which the set $D$ of vertices labelled \texttt{C} induces a clique
of size $\le\lambda$.  The summand for $H^{(D,U)}$ includes the factor
\[
\llbracket\,\hbox{a $\lambda$-colouring of $H[D]$ forces a $\lambda$-colouring of $H$}\rrbracket,
\]
again using the Iverson bracket.  This is polynomial-time computable, so we could reasonably
call it a final base case and say that forced $\lambda$-colouring polynomials have two levels
of recursion.  But these base cases are much less simple than the base cases for other
recursions we have considered (e.g., null graphs, for chromatic polynomials and (eventually) for
Go polynomials).  This time, the computational task for each base case is P-complete \cite{farr94}.

It would be interesting to study this phenomenon of levels of recursion for graph polynomials
more systematically.
Perhaps the notion can be formalised and then related to the logical structure
of the definition of the polynomial.

\section{Graph \textit{polynomials}?}
\label{sec:gph-polys-qn}

We perhaps take it for granted that graph invariants giving counts, or probabilities, of structures
of interest are \textit{polynomials}.  This is not a \textit{necessary} feature of such invariants.
In general, the $\lambda$-Whitney function of a graph \cite{farr04} may have irrational exponents,
though in certain forms this can be avoided when evaluating them along hyperbolae $xy=2^r$ for
$r\in\mathbb{N}$.
(Some related polynomials in knot theory may have negative or fractional exponents, but this is
not a significant exception because in general they can be transformed to polynomials by
appropriate changes of variable.)

The choice of parameter is crucial.  For example, define $\hbox{HomCyc}(G;q)$ to be the
number of homomorphisms from $G$ onto the cycle $C_q$.  Such homomorphisms may be
viewed as $q$-assignments
in which adjacent vertices in $G$ are mapped to ``colours'' (being vertices in $C_q$) that are
adjacent in $C_q$; unlike normal graph colouring (when the homomorphism is to $K_q$),
the two distinct colours used on the endpoints of an edge cannot be completely arbitrary but
are constrained to be neighbouring colours (vertices) in $C_q$.
This is not, in general, a polynomial in $q$.
One way to see this is to note that $\hbox{HomCyc}(K_3;q)=0$ for all even $q$ but is not identically 0.
In other words, in the terminology of de la Harpe and Jaeger \cite{delaharpe-jaeger1995},
the sequence $(C_n:n\in\mathbb{N})$ is not a \textit{strongly polynomial sequence} of graphs.
See \cite{goodall-nesetril-ossonademendez2016} for a study of the deep question of which
sequences of graphs, treated as targets of homomorphisms, give rise to graph polynomials that count
homomorphisms.

Graph colouring requires the colour classes to induce null graphs, where the chromons each
consist of a single vertex.  There has been a lot of work
on generalised colourings where the chromons are less severely restricted.  One of the simplest
relaxations is to bound the sizes of the chromons.
Define $\hbox{mc}(G;s)$ to be the number of 2-assignments of $G$ in which every
chromon has size $\le s$.  We have
\begin{eqnarray*}
\hbox{mc}(G;0)  & = &
\left\{
\begin{array}{cl}
1,  &  \hbox{if $n=0$},   \\
0,  &  \hbox{if $n\ge1$};
\end{array}
\right.   \\
\hbox{mc}(G;1)  & = &
\left\{
\begin{array}{cl}
2^{k(G)},  &  \hbox{if $G$ is bipartite},   \\
0,  &  \hbox{otherwise};
\end{array}
\right.   \\
\hbox{mc}(G;s)  & = &  2^n ~~~\forall s\ge n.
\end{eqnarray*}
So $\hbox{mc}(G;s)$ cannot be a polynomial in $s$.

Define $\hbox{emb}(G;g)$ to be the number of orientable 2-cell 
combinatorial embeddings of $G$ of genus $g$, where embeddings are given by rotation schemes.
This cannot be a polynomial in $g$ because it is positive when $g$ lies between the genus and
maximum genus, inclusive, of $G$, but is zero for all $g$ above the maximum genus.

For graph invariants that \textit{are} polynomials,
it is natural to hope that the theory of polynomials may shed light on the graph polynomials
and, through them, on graphs themselves.  After all, polynomials have a rich mathematical
theory that has been built up over centuries.

For example, Birkhoff's work on the chromatic polynomial, beginning with \cite{birkhoff1912-1913},
was motivated by the thought that its properties, as a polynomial,
might help prove the Four-Colour Conjecture (as it then was).
According to Morse~\cite[p.~386]{morse1946},
``Birkhoff hoped that the theory of chromatic polynomials could be
so developed that methods of analytic function theory could be applied.''

With this motivation, we can ask of any graph polynomial, which aspects of the theory
of polynomials correspond to properties of the underlying graph?

It is common for graph polynomials to be multiplicative over components or even over blocks
(with multiplicativity over blocks being typical for polynomials that depend only on the cycle
matroid of the graph, such as the Tutte polynomial).
This is to be expected for polynomials that count things or give probabilities,
since the lack of interaction between separate components or blocks
typically means we can treat them as contributing independently to counts or probabilities.
Does such a relationship work both ways?  In other words, does multiplicativity \textit{only} occur
over components/blocks?  What, in the graph, is represented by the polynomial's proper factors?
How can we characterise the structure of graphs whose polynomial is irreducible?

In the case of the Tutte polynomial, it was shown by Merino, de Mier and Noy
that the Tutte polynomial of a matroid is irreducible if and only if
the matroid is connected \cite{merino-deMier-noy01}.
So the factors of the polynomial correspond exactly to the components of the matroid,
which means that, for graphs, they correspond to blocks.

But the situation is not so straightforward for many other graph polynomials, even for those
that are specialisations of the Tutte polynomial.

For the chromatic polynomial, it was known to Whitney \cite[\S14]{whitney32d} that the
chromatic polynomial factorises when the graph is \textit{clique-separable}\footnote{but keep in mind
that this term also has a completely different meaning \cite{gavril77}} in the sense that
it has a separating clique, where
a \textit{separating clique} is a clique whose removal increases the number
of components of the graph.
If $G$  is formed from overlap of $H_1$ and $H_2$ in a separating $r$-clique, then Whitney
showed that
\[
P(G;x) = \frac{P(H_1;x)P(H_2;x)}{P(K_r;x)} \,.
\]
Somewhat surprisingly, such \textit{chromatic factorisations} can occur in other cases, too:
in \cite{morgan-farr2009a,morgan-farr2009b}, examples are given of chromatic factorisations
of graphs that are \textit{strongly non-clique-separable} in that they are not chromatically
equivalent to a clique-separable graph.
Some studies of this kind have since been done for other polynomials including the
reliability polynomial \cite{morgan-chen2017} and the stability polynomial \cite{mo-farr-morgan2014}.

Another fundamental topic in the mathematical theory of polynomials is Galois theory.
So it is natural to ask about the relationship between the structure of a graph and the Galois group
of a graph polynomial derived from it.  An initial investigation of this topic for the chromatic
polynomial, including computational results, is reported in \cite{morgan2012}.

The most fundamental property in any mathematical system is \textit{identity}: when are two
objects considered the same?  For a graph polynomial, this is when two graphs
are equivalent in the sense that they have the same polynomial.  This leads to the notion of
certificates of equivalence, which can also be adapted to chromatic factorisation and which we
discuss in the next section.

\section{Certificates}
\label{sec:certs}

The notion of a \textit{certificate} to explain chromatic factorisation and chromatic equivalence
was introduced by Morgan and Farr \cite{morgan-farr2009a,morgan-farr2009b} and developed
further in \cite{morgan2010a,morgan2010b,bukovac-farr-morgan2019}.  The idea has since
been extended to other graph polynomials including
the stability polynomial \cite{mo-farr-morgan2014}, reliability polynomial \cite{morgan-chen2017}
and Tutte polynomial \cite{mo-farr-morgan2014}.
This use of the term ``certificate'' is inspired by its use in complexity theory, e.g., in defining NP,
but we are using the term in a much more specific sense.

Informally, a \textit{certificate} is a sequence of expressions
$E_0, E_1, E_2, \ldots, E_k$ in graphs
where each expression $E_i$, $i>0$, can be obtained from its
predecessor $E_{i-1}$ by applying a relation satisfied by the graph polynomial in question
(e.g., a deletion-contraction relation, or multiplicativity for disjoint unions).
In these expressions, a graph can be regarded as representing its corresponding polynomial.  Replacing each graph by its polynomial, then simplifying the entire
expression, gives a polynomial that can be thought of as the graph polynomial for that expression.

\begin{figure}[h]
    \centering
    \includegraphics[scale=0.4]{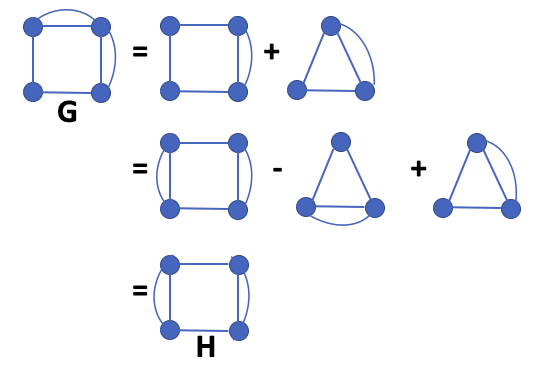}
    \caption{Certificate of Tutte Equivalence. 
 Graph $G$ is Tutte Equivalent to Graph $H$.}
    \label{figTE}
\end{figure}

An example is given in Figure \ref{figTE}.
In this certificate, we use the deletion-contraction relation for the Tutte polynomial in the form 
$$\textbf{[T1]:} ~~~~~~~~ G \longrightarrow G\setminus e + G/e,$$
with each graph standing for its Tutte polynomial.

A simple rearrangement of $T1$ and renaming of the graphs used in the terms gives us
$$\textbf{[T2]:} ~~~~~~~~ G\longrightarrow(G+uv) - G/uv .$$

The resulting certificate has the form: 
\begin{align}
\label{eq:cert-Tutte-equiv}
    G& = G\setminus e + G/e ~~~~~~~~ \text{ (Applying \textbf{T1})}\nonumber\\
    & = (G\setminus e + uv) - G\setminus e /uv + G/e ~~~~~~~~ \text{ (Applying \textbf{T2}})\nonumber \\
    & = G\setminus e + uv ~~~~~~~~ \text{ (Algebraic Step (cancellation))}
\end{align}
where the graph $G\setminus e + uv\cong H$.  If we replace each graph by its Tutte polynomial, each expression in the sequence is equal to the Tutte polynomial of $G$, and hence we have a certificate that shows that graphs $G$ and $H$ are Tutte equivalent. 

Our work so far has focused on certificates of equivalence and factorisation.
We expect that certificates could provide a graph-theoretic approach to studying other algebraic properties
of graph polynomials.

In \cite{morgan2010a}, Morgan introduced the concept of a \textit{schema} or template for
certificates of factorisation and equivalence.
A schema specifies the structure of a certificate, including the relation
to be applied at each step, but without filling in the actual graphs.
So, instead of actual graphs (as in Figure~\ref{figTE}), we just have symbols for them.
In fact, as written --- with symbols $G$, $G\setminus e$, etc.\ --- (\ref{eq:cert-Tutte-equiv})
is really a schema for certificates, and the certificate in Figure~\ref{figTE} is one particular certificate
that belongs to this schema.

In \cite{morgan-farr2009a,morgan-farr2009b}, it was shown that
every graph in a particular infinite family of strongly non-clique-separable graphs
has a chromatic factorisation.  Each certificate
of factorisation in this family used the same sequence of certificate steps,
so the entire infinite family of certificates could be described by a single schema.

If two graphs have the same multiset of blocks, then they are Tutte equivalent, since the Tutte
polynomial is multiplicative over blocks.  We can capture this multiplicativity in certificate steps
that allow a graph to be replaced by a formal product of its blocks and vice versa.
\begin{description}
    \item[\textbf{[T3]}:] $G\longrightarrow B_1  B_2  \cdots B_k$ where the $B_i$ are the blocks of $G$,
    \item[\textbf{[T4]}:]  $B_1  B_2  \cdots B_k\longrightarrow G$ where $G$ is a graph with blocks $B_i$, $1\leq i\leq k$,
\end{description}
This enables us to write the following simple schema for certificates of Tutte equivalence for pairs of graphs with the same blocks.
\begin{align}
    G & = \prod_{i=1}^{k} B_i ~~~~~~~~ \text{ (Applying \textbf{T3})}\nonumber\\
    & = H ~~~~~~~~ \text{ (Applying \textbf{T4}).}
\end{align}
Effectively, this schema first `unglues' blocks then 'glues' them back together to produce graph $H$.  This certificate schema works for all pairs of graphs that have the same blocks and may be
regarded as a representation of the set of all such pairs.

Schemas 1 and 2 in \cite{bukovac-farr-morgan2019} give two of the shortest certificates for pairs of chromatically equivalent graphs, $G$ and $H$.  In both these schemas, graph $H$ can be obtained by removing an edge from graph $G$ and then adding a different edge.

Schema 2 relates pairs of chromatically equivalent graphs $G$ and $H$ with $G\setminus e\cong H\setminus f$ and $G/e\cong H/f$, where $e\in E(G)$ and $f\in E(H)$.  Applying two deletion-contraction steps, we have: 
\begin{align}
    G&= G\setminus e - G/e \nonumber\\
    & = (G\setminus e) + f\label{stepRevDC}
\end{align}
where $f\not\in E(G)$ and $G\setminus e+f\cong H$.  The second step, (\ref{stepRevDC}), is obtained by rearranging the usual deletion-contraction relation.  

Schema 1 is similar to Schema 2, but uses the addition-identification relation.  Here $(G+e)\setminus f\cong H$.  Applying two addition-identification steps, we have:
\begin{align}
    G& = G + e + G/e \nonumber \\
    & = (G+e)\setminus f \label{stepRevAI}
\end{align}
where $e\not\in E(G)$ and $f\in E(G)$. 

More sophisticated certificates of equivalence are available.  In \cite{bukovac-farr-morgan2019},
shortest certificates of chromatic equivalence are given for all pairs of chromatically equivalent graphs of order at most 7.  These corresponded to 15 different schemas.  It should be noted that a shortest certificate of equivalence may not be unique.  In \cite{morgan2010b}, infinitely many pairs of chromatically equivalent non-isomorphic graphs are constructed along with their
certificates of equivalence.

The length of certificates has implications for the complexity of testing equivalence with respect
to these polynomials (chromatic equivalence, Tutte equivalence, etc.).
A short certificate of equivalence, once obtained, gives a means of verifying that two graphs
have the same polynomial without computing the polynomial.  For example, if the length of
certificates of chromatic equivalence is polynomially bounded, then the problem of testing
chromatic equivalence belongs to NP \cite{bukovac-farr-morgan2019}.
But at present we only have very loose,
exponential upper bounds on certificate length.  We discuss some implications of certificate length
for the computational complexity of chromatic equivalence, factorisation and uniqueness
in \cite{bukovac-farr-morgan2019} and Tutte equivalence in \cite{mo-morgan-farrXX}.

Appropriate versions of certificates of equivalence and factorisation should be applicable
to many other graph polynomials.
We would like to see a rigorous theory of certificates of (at least) equivalence
and factorisation for a broad class of graph polynomials with reduction relations.
The commutativity of the operations of deleting/contracting different
edges in a graph may be regarded as an instance of the Church-Rosser property
from the theory of rewriting systems, as observed by Yetter \cite{yetter1990} and Makowsky \cite{makowsky08}; see also \cite[result 9m, p.~72]{biggs1993}.
It seems to us that the theory of rewriting systems could shed more light
on the kind of certificates we have considered here.   \\

\noindent\textbf{Acknowledgements}

We thank Andrew Goodall, J\'anos Makowsky, Steven Noble and the referee for their comments.
Through this \textit{Festschrift} contribution,
it is a pleasure to acknowledge J\'anos's far-reaching contributions to the
study of graph polynomials, through his mathematics and also through his generosity
as a colleague, in sharing ideas, organising meetings and supporting the
work of others.

\end{document}